\documentclass[11pt]{amsart}

\usepackage[papersize={8.5in,11in},margin=1in]{geometry}

\usepackage{latexsym,eucal,amsmath,amssymb}
\usepackage{fullpage}
\usepackage{color}

\newtheorem{Theorem}{Theorem}[section]
\newtheorem{Lemma}[Theorem]{Lemma}
\newtheorem{Proposition}[Theorem]{Proposition}
\newtheorem{Corollary}[Theorem]{Corollary}
\newtheorem{Conjecture}[Theorem]{Conjecture}
\newtheorem{Question}[Theorem]{Question}
\newtheorem{Definition}[Theorem]{Definition}

\newenvironment{Remark}{\begin{trivlist} \item[] {\bf Remark.}}{\end{trivlist}}
\newenvironment{Proof}{\begin{trivlist} \item[] {\bf Proof.}}{\hfill $\Box$\end{trivlist}}

\newcommand{\sB}{\EuScript{B}}

\renewcommand{\geq}{\geqslant}
\renewcommand{\leq}{\leqslant}

\begin{document}
\title[Colourings of path systems]{Colourings of path systems}

\author{Iren Darijani and  David A. Pike\\
Department of Mathematics and Statistics\\
Memorial University of Newfoundland\\
St.\ John's, NL, Canada. A1C 5S7}

\maketitle

\begin{abstract}
A $P_m$ path in a graph is a path on $m$ vertices.  A $P_m$ system of order $n>1$ is a partition of the edges of the complete graph $K_n$ into $P_m$ paths.  A $P_m$ system is said to be $k$-colourable if the vertex set of $K_n$ can be partitioned into $k$ sets called colour classes such that no path in the system is monochromatic. The system is $k$-chromatic if it is $k$-colourable but is not $(k-1)$-colourable. If every $k$-colouring of a $P_m$ system can be obtained from some $k$-colouring $\phi$ by a permutation of the colours, we say that the system is uniquely $k$-colourable. In this paper, we first observe that there exists a $k$-chromatic $P_m$ system for any $k\geq 2$ and $m\geq 4$ where $m$ is even. Next, we prove that there exists an equitably 2-chromatic $P_4$ system of order $n$ for each admissible order $n$. We then show that for all $k\geq 3$, there exists a $k$-chromatic $P_4$ system of order $n$ for all sufficiently large admissible $n$. Finally, we show that there exists a uniquely 2-chromatic $P_4$ system of order $n$ for each admissible $n \geq 109$.
\end{abstract}

\vspace*{\baselineskip}

\vspace*{\baselineskip}

\section{Introduction}
A {\em $G$-decomposition} of a graph $H$ is a pair $(V,\mathcal{B})$ where $V$ is the set of vertices of $H$ and $\mathcal{B}$ is a set of subgraphs of $H$, each isomorphic to $G$, whose edge sets partition the edge set of $H$. A {\em $G$-design} of order $n$ is a $G$-decomposition of the complete graph $K_n$ on $n$ vertices.  A $G$-design $(V,\mathcal{B})$  is said to be {\em weakly $k$-colourable} if its vertex set can be partitioned into $k$ sets (called colour classes) such that no subgraph belonging to $\sB$ is monochromatic. The  $G$-design is {\em $k$-chromatic} if it is $k$-colourable but is not $(k-1)$-colourable. If a $G$-design is $k$-chromatic, we say that its chromatic number is $k$. A colouring of a  $G$-design is said to be {\em equitable} if the cardinalities of the colour classes differ by at most one. It is {\em strongly equitable} if the colour classes are of the same size. A $G$-design will be called {\em (strongly) equitably $k$-chromatic} if it is $k$-chromatic and admits a (strongly) equitable $k$-colouring. If every $k$-colouring of a $G$-design $(V,\mathcal B)$ can be obtained from some $k$-colouring $\phi$ by a permutation of the colours, we say that $(V,\mathcal B)$ has a {\em unique} $k$-colouring, or that $(V,\mathcal B)$ is {\em uniquely $k$-colourable}.

Weak colourings of $G$-designs have been extensively studied for $K_3$-designs, i.e., Steiner triple systems. In particular, de Brandes, Phelps, and R\"{o}dl proved that for any $k\geq 3$, there exists some integer $v_k$ such that for all admissible $v\geq v_k$, there is a $k$-chromatic Steiner triple system of order $v$~\cite{de Brandes}.
In 2008, Burgess and Pike showed that for all $k\geq 2$ and even $m\geq 4$ there exists a $k$-chromatic $m$-cycle system~\cite{Burgess}.
In 2010, Horsley and Pike showed that for all $k\geq 2$ and $m\geq 3$ with $(k,m)\neq (2,3)$, there exist $k$-chromatic $m$-cycle systems of all admissible orders greater than or equal to some integer $n_{k,m}$~\cite{HorsleyPike2010}.  Horsley and Pike subsequently considered balanced incomplete block designs and showed in 2014 that the obvious necessary conditions for the existence of a $(v,m,\lambda)$-BIBD are asymptotically sufficient for the existence of a $k$-chromatic $(v,m,\lambda)$-BIBD except when $(k,m)=(2,3)$~\cite{HorsleyPike2014}.

More recently, Darijani and Pike~\cite{Dar} showed that for any integer $k\geq 2$, there exists a $k$-chromatic 3-star system of order $n$ for all sufficiently large admissible $n$. They also generalized this result for $e$-star systems for any $e\geq 3$. They showed that for all $k\geq 2$ and $e\geq 3$, there exists a $k$-chromatic $e$-star system of order $n$ for all sufficiently large $n$ such that $n\equiv 0,1$ (mod $2e$). An $e$-star is a complete bipartite graph $K_{1,e}$. An $e$-star system of order $n>1$ is a partition of the edges of the complete graph $K_n$ into $e$-stars. They also proved that for all $k\geq 2$ and $e\geq 3$, there exists a uniquely $k$-chromatic $e$-star system of order $n$ for all sufficiently large $n$ such that  $n\equiv 0,1$ (mod $2e$). The idea of uniquely colourable designs has already arisen in the context of Steiner triple systems. In 2003, Forbes showed that for every admissible $v\geq 25$, there exists a 3-balanced Steiner triple system with a unique 3-colouring and also a Steiner triple system  which has a unique nonequitable 3-colouring (where a Steiner triple system is said to be {\em 3-balanced} if every 3-colouring of it is equitable)~\cite{Forbes}.

In this paper, we investigate $k$-colourings of path systems.
For any integer $m\geq 2$, the path with $m$ vertices $\{a_1,\ldots,a_m\}$ is denoted by $P_m=(a_1,\ldots,a_m)$. When $G$ is a $P_m$ path, a $G$-design of order $n$ is called a {\em $P_m$ system} of order $n$. The necessary and sufficient conditions for the existence of a $P_m$ system of  order $n$ are that $n=1$ or $n\geq m$ and $n(n-1)\equiv 0$ (mod $2m-2$)~\cite{Tarsi}; for fixed $m$ such values of $n$ will be called {\em admissible}.  In Section~2, we construct some small systems that are used as ingredients to construct larger systems. In Section~3, we first observe that there exists a $k$-chromatic $P_m$ system for any $k\geq 2$ and $m\geq 4$ where $m$ is even. Next, we prove that there exists an equitably 2-chromatic $P_4$ system for each admissible order $n$ (i.e., for $n \equiv 0,1$ (mod $3$) and $n \geq 4$). We then show that for any  integer $k \geq 3$, there exists some integer $n_k$ such that for all admissible $n \geq n_k$, there exists a $k$-chromatic $P_4$ system of order $n$. In Section~4, we study unique 2-colourings of $P_4$ systems. We first show that there exists a strongly equitably 2-chromatic $P_4$ system of order 28 for which two particular vertices are forced to have different colours.
We conclude by showing that there exists a uniquely 2-chromatic $P_4$ system of order $n$ for each $n \equiv 0,1$ (mod 3), $n \geq 109$.

\section{Some small explicit systems}
In this section, we first construct (strongly) equitably 2-chromatic $P_4$ systems of order 4, 6, and 7. We then present $P_4$-decompositions of several other graphs of small order. These decompositions are used as ingredients to construct larger systems in the next section.
\begin{Lemma}\label{lem0.1}
There exist strongly equitably 2-chromatic $P_4$ systems of order 4 and 6, and an equitably 2-chromatic $P_4$ system of order 7.
\end{Lemma}
\begin{Proof}
We first construct a strongly equitably 2-chromatic $P_4$ system of order 4. Let $\hat{V}=\{v_1,v_2,v_3,v_4\}$, $\hat{\mathcal B}=\big \{ (v_3,v_1,v_2,v_4)$, $(v_1,v_4,v_3,v_2)\big\}$, $\hat{R}=\{v_1,v_3\}$, and $\hat{Y}=\{v_2,v_4\}$. Then $(\hat{V},\hat{\mathcal B})$ is a $P_4(4)$ for which $\hat{R}$ and $\hat{Y}$ constitute a strongly equitable 2-colouring.

Next, we construct a strongly equitably 2-chromatic $P_4$ system of order 6. Let $\tilde{V}=\{v_1,v_2,v_3,v_4,v_5,v_6\}$, $\tilde{\mathcal B}=\big\{(v_6,v_1,v_2,v_5)$, $(v_6,v_2,v_3,v_5)$, $(v_2,v_4,v_1,v_3)$, $(v_4,v_3,v_6,v_5)$, $(v_1,v_5,v_4,v_6)\big\}$, $\tilde{R}=\{v_1,v_3,v_5\}$, and $\tilde{Y}=\{v_2,v_4,v_6\}$. Then $(\tilde{V},\tilde{\mathcal B})$ is a $P_4(6)$ for which $\tilde{R}$ and $\tilde{Y}$ constitute a strongly equitable 2-colouring.

Finally, we construct an equitably 2-chromatic $P_4$ system of order 7. Let $V'=\{v_1,v_2,v_3,v_4,v_5,v_6,v_7\}$, $\mathcal B'=\big\{ (v_2,v_7,v_4,v_3)$, $(v_7,v_3,v_6,v_5)$, $(v_1,v_7,v_6,v_2)$, $(v_2,v_4,v_1,v_3)$, $(v_1,v_5,v_4,v_6)$, $(v_2,v_3,v_5,v_7)$, $(v_6,v_1,v_2,v_5)\big\}$, $R'=\{v_1,v_3,v_5,v_7\}$, and $Y'=\{v_2,v_4,v_6\}$. Then $(V',\mathcal B')$ is a $P_4(7)$ for which $R'$ and $Y'$ constitute an equitable 2-colouring.
\end{Proof}

\begin{Lemma}\label{lem.a}
There exists a strongly equitably 2-chromatic decomposition of $K_{3,3}$ into $P_4$ paths.
\end{Lemma}
\begin{Proof}
Let $V=\{v_1,v_2,v_3\}\cup \{w_1,w_2,w_3\}$ be the set of vertices, $\mathcal B=\big\{(v_2,w_1,v_1,w_3)$, $(v_1,w_2,v_3,w_1)$, $(v_3,w_3,v_2,w_2)\big\}$, $R=\{v_1,v_3,w_2\}$, and $Y=\{v_2,w_1,w_3\}$. Then $\mathcal B$ is a decomposition of $K_{3,3}$ into $P_4$ paths for which $R$ and $Y$ constitute a strongly equitable 2-colouring.
\end{Proof}

\begin{Lemma}\label{lem.b}
There exists an equitably 2-chromatic decomposition of $K_{4,3}$ into $P_4$ paths.
\end{Lemma}
\begin{Proof}
Let $V=\{v_1,v_2,v_3,v_4\}\cup \{w_1,w_2,w_3\}$ be the set of vertices, $\mathcal B=\big\{(v_4,w_3,v_1,w_1)$, $(v_2,w_1,v_4,w_2)$, $(v_1,w_2,v_3,w_1)$, $(v_3,w_3,v_2,w_2) \big\}$, $R=\{v_1,v_3,w_1,w_3\}$, and $Y=\{v_2,v_4,w_2\}$. Then $\mathcal B$ is a decomposition of $K_{4,3}$ into $P_4$ paths for which $R$ and $Y$ constitute an equitable 2-colouring.
\end{Proof}

\begin{Lemma}\label{lem.c}
There exists a strongly equitably 2-chromatic decomposition of $K_{7,3}$ into $P_4$ paths.
\end{Lemma}
\begin{Proof}
Let $V=\{v_1,v_2,v_3,v_4,v_5,v_6,v_7\}\cup \{w_1,w_2,w_3\}$ be the set of vertices, $\mathcal B=\big\{(v_1,w_1,v_7,w_3)$, $(v_7,w_2,v_2,w_1)$, $(v_2,w_3,v_3,w_1)$, $(v_3,w_2,v_4,w_1)$, $(v_4,w_3,v_5,w_1)$,  $(v_5,w_2,v_6,w_1)$, $(v_6,w_3,v_1,w_2) \big\}$, $R=\{v_1,v_3,v_5,v_7,w_2\}$, and $Y=\{v_2,v_4,v_6,w_1,w_3\}$. Then $\mathcal B$ is a decomposition of $K_{4,3}$ into $P_4$ paths for which $R$ and $Y$ constitute a strongly equitable 2-colouring.
\end{Proof}

\begin{Lemma}\label{lem.d}
There exists a strongly equitably 2-chromatic decomposition of $K_{6,2}$ into $P_4$ paths.
\end{Lemma}
\begin{Proof}
Let $V=\{ v_1,v_2,v_3,v_4,v_5,v_6\}\cup\{w_1,w_2\}$ be the set of vertices, $\mathcal B=\big\{(w_1,v_1,w_2,v_2)$, $(v_2,w_1,v_3,w_2)$, $(w_1,v_4, w_2,v_5)$, $(v_5,w_1,v_6,w_2)\big\}$, $R=\{v_1,v_3,v_5,w_1\}$, and $Y=\{v_2,v_4,v_6,w_2\}$. Then $\mathcal B$ is a decomposition of $K_{6,2}$ into $P_4$ paths for which $R$ and $Y$ constitute a strongly equitable 2-colouring.
\end{Proof}

\begin{Lemma}\label{lem.e}
There exists an equitably 2-chromatic decomposition of $K_{6,3}$ into $P_4$ paths.
\end{Lemma}
\begin{Proof}
Let $V=\{v_1,v_2,v_3,v_4,v_5,v_6\}\cup\{w_1,w_2,w_3\}$ be the set of vertices, $\mathcal B=\big\{(v_1,w_1,v_2,w_2)$, $(v_2,w_3,v_3,w_1)$, $(v_3,w_2,v_4,w_1),$ $(v_4,w_3,v_5,w_1)$, $(v_5,w_2,v_6,w_1)$, $(v_6,w_3,v_1,w_2)\big\}$, $R=\{v_1,v_3,v_5,w_1,w_3\}$, and $Y=\{v_2,v_4,v_6,w_2\}$. Then $\mathcal B$ is a decomposition of $K_{6,3}$ into $P_4$ paths for which $R$ and $Y$ constitute an equitable 2-colouring.
\end{Proof}

\begin{Lemma}\label{lem.f}
There exists a strongly equitably 2-chromatic decomposition of $K_{6,4}$ into $P_4$ paths.
\end{Lemma}
\begin{Proof}
Let $V=\{v_1,v_2,v_3,v_4,v_5,v_6\}\cup\{w_1,w_2,w_3,w_4\}$ be the set of vertices, $\mathcal B=\big\{(v_1,w_1,v_2,w_2)$, $(v_1,w_3,v_2,w_4)$, $(v_1,w_2,v_3,w_1)$, $(v_1,w_4,v_3,w_3)$, $(v_4,w_1,v_5,w_2)$, $(v_4,w_3,v_5,w_4)$, $(v_4,w_2,v_6,w_1)$, $(v_4,w_4,v_6,w_3)\big\}$, $R=\{v_1,v_3,v_5,w_1,w_3\}$, and $Y=\{v_2,v_4,v_6,w_2,w_4\}$. Then $\mathcal B$ is a decomposition of $K_{6,4}$ into $P_4$ paths for which $R$ and $Y$ constitute a strongly equitable 2-colouring.
\end{Proof}

\begin{Lemma}\label{lem.g}
There exists an equitably 2-chromatic decomposition of $K_{6,5}$ into $P_4$ paths.
\end{Lemma}
\begin{Proof}
Let $V=\{v_1,v_2,v_3,v_4,v_5,v_6\}\cup\{w_1,w_2,w_3,w_4,w_5\}$ be the set of vertices, $\mathcal B=\big\{(v_1,w_1,v_2,w_2)$, $(v_2,w_3,v_3,w_1)$, $(v_3,w_2,v_4,w_1)$, $(v_4,w_3,v_5,w_1)$, $(v_5,w_2,v_6,w_1)$, $(v_6,w_3,v_1,w_2)$, $(v_2,w_4,v_3,w_5)$, $(v_2,w_5,v_4,w_4)$, $(v_5,w_4,v_6,w_5)$, $(v_5,w_5,v_1,w_4)\big\}$, $R=\{v_1$, $v_3,v_5,w_1,w_3,w_5\}$, and $Y=\{v_2,v_4,v_6,w_2,w_4\}$. Then $\mathcal B$ is a decomposition of $K_{6,5}$ into $P_4$ paths for which $R$ and $Y$ constitute an equitable 2-colouring.
\end{Proof}

\begin{Lemma}\label{lem.h}
There exists a strongly equitably 2-chromatic decomposition of $K_{6,6}$ into $P_4$ paths.
\end{Lemma}
\begin{Proof}
Let $V=\{v_1,v_2,v_3,v_4,v_5,v_6\}\cup\{w_1,w_2,w_3,w_4,w_5,w_6\}$ be the set of vertices, $\mathcal B=\big\{(v_1,w_1,v_2,w_2)$, $(v_2,w_3,v_3,w_1)$, $(v_3,w_2,v_4,w_1)$, $(v_4,w_3,v_5,w_1)$, $(v_5,w_2,v_6,w_1)$, $(v_6,w_3,v_1,w_2)$, $(v_1,w_4,v_2,w_5)$, $(v_2,w_6,v_3,w_4)$, $(v_3,w_5,v_4,w_4)$, $(v_4,w_6,v_5,w_4)$, $(v_5,w_5,v_6,w_4)$, $(v_6,w_6,v_1,w_5)\big\}$, $R=\{v_1,v_3,v_5,w_1,w_3,w_5\}$, and $Y=\{v_2,v_4,v_6,w_2,w_4,w_6\}$. Then $\mathcal B$ is a decomposition of $K_{6,6}$ into $P_4$ paths for which $R$ and $Y$ constitute a strongly equitable 2-colouring.
\end{Proof}

\begin{Lemma}\label{lem.i}
There exists an equitably 2-chromatic decomposition of $K_{6,7}$ into $P_4$ paths.
\end{Lemma}
\begin{Proof}
Let $V=\{v_1,v_2,v_3,v_4,v_5,v_6\}\cup\{w_1,w_2,w_3,w_4,w_5,w_6,w_7\}$ be the set of vertices, $\mathcal B=\big\{(v_1,w_1,v_2,w_2)$, $(v_2,w_3,v_3,w_1)$, $(v_3,w_2,v_4,w_1)$, $(v_4,w_3,v_5,w_1)$, $(v_5,w_2,v_6,w_1)$, $(v_6,w_3,v_1,w_2)$, $(v_1,w_4,v_2,w_5)$, $(v_1,w_6,v_2,w_7)$, $(v_1,w_5,v_3,w_4)$, $(v_1,w_7,v_3,w_6)$, $(v_4,w_4,v_5,w_5)$, $(v_4,w_6,v_5,w_7)$, $(v_4,w_5,v_6,w_4)$, $(v_4,w_7,v_6,w_6)\big\}$, $R=\{v_1,v_3,v_5,w_1,w_3,w_5,w_7\}$, and $Y=\{v_2,v_4,v_6,w_2,w_4,w_6\}$. Then $\mathcal B$ is a decomposition of $K_{6,7}$ into $P_4$ paths for which $R$ and $Y$ constitute an equitable 2-colouring.
\end{Proof}

\begin{Lemma}\label{lem6}
Let $(V_1,V_2)$ be a bipartition of the set of vertices of the complete bipartite graph $K_{6,3}$ where $|V_1|=6$ and $|V_2|=3$. Let $E_1$ be the set of edges of $K_{6,3}$ and $E_2$ be the set of edges of the complete graph on $V_2$. Then the graph on vertex set $V_1 \cup V_2$ with edge set $E_1\cup E_2$ has a decomposition into $P_4$ paths which is equitably 2-chromatic.
\end{Lemma}
\begin{Proof}
Let $V_1=\{v_1,v_2,v_3,v_4,v_5,v_6\}$, $V_2=\{w_1,w_2,w_3\}$, and $\mathcal B=\big\{(v_1,w_1,w_2,w_3)$, $(w_2,v_2,w_1,w_3),$ $(v_2,w_3,v_3,w_1)$ $,(v_3,w_2,v_4,w_1),$ $(v_4,w_3,v_5,w_1),$ $(v_5,w_2,v_6,w_1),$ $(v_6,w_3,v_1,w_2)\big\}$, $R=\{v_1,v_3,v_5,w_1,w_3\}$, and $Y=\{v_2,v_4,v_6,w_2\}$. Then $\mathcal B$ is a decomposition of $E_1\cup E_2$ into $P_4$ paths for which $R$ and $Y$ constitute an equitable 2-colouring.
\end{Proof}

\begin{Lemma}\label{lem7}
Let $(V_1,V_2)$ be a bipartition of the set of vertices of the complete bipartite graph $K_{7,2}$ where $|V_1|=7$ and $|V_2|=2$. Let $E_1$ be the set of edges of $K_{7,2}$ and $E_2$ be the set of edges of the complete graph on $V_2$. Then the graph on vertex set $V_1 \cup V_2$ with edge set $E_1\cup E_2$ has a decomposition into $P_4$ paths which is equitably 2-chromatic.
\end{Lemma}
\begin{Proof}
Let $V_1=\{v_1,v_2,v_3,v_4,v_5,v_6,v_7\}$, $V_2=\{w_1,w_2\}$, and $\mathcal B=\big\{(v_2,w_1,v_3,w_2),$ $(v_2,w_2,w_1,v_4),$ $(v_4,w_2,v_5,w_1),$ $(v_6,w_1,v_7,w_2),(v_6,w_2,v_1,w_1)\big\}$, $R=\{v_1,v_3,v_5,v_7,w_1\}$, and $Y=\{v_2,v_4,v_6,w_2\}$. Then $\mathcal B$ is a decomposition of $E_1\cup E_2$ into $P_4$ paths for which $R$ and $Y$ constitute an equitable 2-colouring.
\end{Proof}

\begin{Lemma}\label{lem8}
Let $(V_1,V_2)$ be a bipartition of the set of vertices of the complete bipartite graph $K_{7,3}$ where $|V_1|=7$ and $|V_2|=3$. Let $E_1$ be the set of edges of $K_{7,3}$ and $E_2$ be the set of edges of the complete graph on $V_2$. Then the graph on vertex set $V_1 \cup V_2$ with edge set $E_1\cup E_2$ has a decomposition into $P_4$ paths which is strongly equitably 2-chromatic.
\end{Lemma}
\begin{Proof}
Let $V_1=\{v_1,v_2,v_3,v_4,v_5,v_6,v_7\}$, $V_2=\{w_1,w_2,w_3\}$, and $\mathcal B=\big\{(v_4,w_1,v_5,w_3)$, $(v_5,w_2,v_4,w_3),$ $(v_6,w_1,v_7,w_3),$ $(v_7,w_2,v_6,w_3),(v_1,w_3,v_3,w_1)$, $(v_2,w_1,v_1,w_2)$, $(v_3,w_2,w_3,v_2)$, $(v_2,w_2,w_1,w_3)\big\}$, $R=\{v_1,v_3,v_5,v_7,w_1\}$, and $Y=\{v_2,v_4,v_6,w_2,w_3\}$. Then $\mathcal B$ is a decomposition of $E_1\cup E_2$ into $P_4$ paths for which $R$ and $Y$ constitute a strongly equitable 2-colouring.
\end{Proof}

\begin{Lemma}\label{lem13}
Let $(V_1,V_2)$ be a bipartition of the set of vertices of the complete bipartite graph $K_{4,2}$ where $|V_1|=4$ and $|V_2|=2$. Let $E_1$ be the set of edges of $K_{4,2}$ and $E_2$ be the set of edges of the complete graph on $V_2$. Then the graph on vertex set $V_1 \cup V_2$ with edge set $E_1\cup E_2$ has a decomposition into $P_4$ paths which is strongly equitably 2-chromatic.\end{Lemma}
\begin{Proof}
Let $V_1=\{v_1,v_2,v_3,v_4\}$, $V_2=\{w_1,w_2\}$, and $\mathcal B=\big\{(v_1,w_1,w_2,v_4),(v_4,w_1,v_2,w_2)$, $(v_1,w_2,v_3,w_1)\big\}$, $R=\{v_1,v_3,w_1\}$, and $Y=\{v_2,v_4,w_2\}$. Then $\mathcal B$ is a decomposition of $E_1\cup E_2$ into $P_4$ paths for which $R$ and $Y$ constitute a strongly equitable 2-colouring.
\end{Proof}

\begin{Lemma}\label{lem14}
Let $(V_1,V_2)$ be a bipartition of the set of vertices of the complete bipartite graph $K_{4,5}$ where $|V_1|=4$ and $|V_2|=5$. Let $E_1$ be the set of edges of $K_{4,5}$ and $E_2$ be the set of edges of the complete graph on $V_2$. Then the graph on vertex set $V_1 \cup V_2$ with edge set $E_1\cup E_2$ has a decomposition into $P_4$ paths which is equitably 2-chromatic.\end{Lemma}
\begin{Proof}
Let $V_1=\{v_1,v_2,v_3,v_4\}$, $V_2=\{w_1,w_2,w_3,w_4,w_5\}$, and $\mathcal B=\big\{(v_4,w_3,v_1,w_1)$, $(v_2,w_1,v_4,w_2)$, $(v_1,w_2,v_3,w_1)$, $(v_3,w_3,v_2,w_2)$, $(v_1,w_4,w_5,v_4)$, $(v_4,w_4,v_2,w_5)$, $(v_1,w_5,v_3,w_4)$, $(w_4,w_1,w_2,w_3)$, $(w_4,w_2,w_5,w_3)$, $(w_4,w_3,w_1,w_5)\big\}$, $R=\{v_1,v_3,w_1,w_3,w_5\}$, and $Y=\{v_2,v_4,w_2,w_4\}$. Then $\mathcal B$ is a decomposition of $E_1\cup E_2$ into $P_4$ paths for which $R$ and $Y$ constitute an equitable 2-colouring.
\end{Proof}

\section{$k$-colourings of path systems}
In this section, we first observe that there exists a $k$-chromatic $P_m$ system for any $k\geq 2$ and $m\geq 4$ where $m$ is even. We then concentrate on $m=4$ and prove that there exists an equitably 2-chromatic $P_4$ system for each admissible order $n$. We finish this section by showing that for any  integer $k \geq 3$, there exists some integer $n_k$ such that for all admissible $n \geq n_k$, there exists a $k$-chromatic $P_4$ system of order $n$.

\subsection{$\mathbf k$-chromatic $\mathbf {P_m}$ systems for $\mathbf m$ even}
In this subsection, we prove that there exists a $k$-chromatic $P_m$ system for any $k\geq 2$ and $m\geq 4$ where $m$ is even. Recall that for any $m\geq 2$, necessary and sufficient conditions for the existence of a $P_m$ system of order $n$ are given in the following theorem from~\cite{Tarsi}.
\begin{Theorem}\cite{Tarsi}\label{Tarsi1}
Let $m\geq 2$. There exists a $P_m$ system of order $n$ if and only if
\begin{itemize}
\item $n=1$ or $n \geq m$; and
\item $n(n-1)\equiv 0$ (mod $2m-2$).
\end{itemize}
\end{Theorem}

To prove the main result of this subsection, we begin with the following theorem from~\cite{HorsleyPike2014}.

\begin{Theorem}\cite{HorsleyPike2014}\label{thm-HP2014}
Let $k$, $n$ and $\lambda$ be positive integers such that $k\geq 2$, $n\geq 3$ and $(k,n)\neq (2,3)$. Then there is an integer $N(k,n,\lambda)$ such that there exists a $k$-chromatic BIBD$(v,n,\lambda)$ for all admissible integers $v\geq N(k,n,\lambda)$.
\end{Theorem}
We now show that any complete graph $K_n$, where $n$ is even, has a decomposition into Hamilton paths.
\begin{Lemma}\label{lem_hamilton}
For any even $n$, the complete graph $K_n$ has a decomposition into Hamilton paths.
\end{Lemma}
\begin{Proof}
The complete graph $K_{n+1}$ has a Hamilton cycle decomposition as demonstrated in the 19th century by Walecki~\cite{Lucas}. Remove one vertex from $K_{n+1}$ to get a decomposition of $K_{n}$ into Hamilton paths. This lemma also follows from Theorem~\ref{Tarsi1}.
\end{Proof}
We then obtain the main result of this subsection.
\begin{Corollary}\label{cor1}
For any $m\geq 4$, $m$ even, and $k\geq 2$, there exists a $k$-chromatic $P_m$ system.
\end{Corollary}
\begin{Proof}
By Theorem~\ref{thm-HP2014}, for any $k\geq 2$, $t\geq 2$, and any sufficiently large admissible integer $v$, there exists a weakly $k$-chromatic BIBD$(v,2t,1)$. By Lemma~\ref{lem_hamilton}, the complete graph $K_{2t}$ has a decomposition into Hamilton paths. Let $m=2t$. Therefore, there exists a $k$-chromatic $P_m$ system for any $k\geq 2$ and even $m\geq 4$.
\end{Proof}

Note that a BIBD$(v,4,1)$ exists if and only if $v \equiv 1,4$ (mod 12). Therefore, by Theorem~\ref{thm-HP2014}, for any $k\geq 2$, there exist $k$-chromatic $P_4$ systems of order $v \equiv 1,4$ (mod 12) for all sufficiently large $v$. In the next two subsections, we show that for any $k\geq 2$, there exist $k$-chromatic $P_4$ systems for all sufficiently large admissible $v$.

\subsection{2-chromatic $\mathbf{P_4}$ systems}
In this subsection, we prove that for each admissible order $n$, there exists an equitably 2-chromatic $P_4$ system of order $n$. Recall that a $P_4$ system of order $n$  exists if and only if $n\equiv 0,1,3,4$ (mod $6$) by Theorem~\ref{Tarsi1}. To prove the main result of this subsection, we initially establish the following two lemmas.
\begin{Lemma}\label{lem0.3}
Let $t\geq 1$. There exist a strongly equitably 2-chromatic $P_4$ system of order $6t$ and an equitably 2-chromatic $P_4$ system of order $6t+1$.
 \end{Lemma}
\begin{Proof}
For $t=1$ refer to Lemma~\ref{lem0.1}. For each $t>1$, we first construct an equitably 2-chromatic $P_4$ system $P_4(6t)$, $(V,\mathcal B)$. Let $V=\{1,\ldots,6t\}$ be the set of points. Partition the set of points $V$ into $t$ subsets $A_1=\{1,2,3,4,5,6\},\ldots$, $A_t=\{6t-5,6t-4,6t-3,6t-2,6t-1,6t\}$.

Note that by Lemma~\ref{lem0.1} and Lemma~\ref{lem.h}, there exist an equitably 2-chromatic $P_4$ system of order six and a decomposition of $K_{6,6}$ into $P_4$ paths which is equitably 2-chromatic. So, for each $1\leq i\leq t$ and $1\leq k<\ell \leq t$, let $(A_i,\mathcal A_i)$ be an equitably 2-chromatic $P_4$ system of order six similar to that constructed from Lemma~\ref{lem0.1}, and let $\mathcal A_{k,\ell}$ be an equitably 2-chromatic decomposition of the edges between $A_k$ and $A_{\ell}$ into $P_4$ paths which is constructed in a similar manner to Lemma~\ref{lem.h}.

Let $\mathcal B=(\bigcup_{i=1}^{t}\mathcal A_i)\cup \big(\bigcup_{k=1}^{t-1}(\bigcup_{\ell=k+1}^{t}\mathcal A_{k,\ell})\big)$, $R=\{1,3,\ldots,6t-1\}$, and $Y=\{2,4,\ldots,6t\}$. Then $(V,\mathcal B)$ is a $P_4(6t)$ for which $R$ and $Y$ constitute a strongly equitable 2-colouring.

Next, we construct an equitably 2-chromatic $P_4$ system $P_4(6t+1)$, $(\hat{V},\hat{\mathcal B})$. Let $\hat{V}=\{1,\ldots,6t+1\}$ be the set of points. Partition the set of points $\hat{V}$ into $t$ subsets $\hat{A_1}=\{1,2,3,4,5,6,7\},\hat{A_2}=\{8,9,10,11,12,13\}, \ldots,\hat{A_t}=\{6t-4,6t-3,6t-2,6t-1,6t,6t+1\}$.

Note that by Lemma~\ref{lem0.1}, Lemma~\ref{lem.h}, and Lemma~\ref{lem.i} there exist equitably 2-chromatic $P_4$ systems of order six and seven and decompositions of $K_{6,6}$ and $K_{6,7}$ into $P_4$ paths which are equitably 2-chromatic. So, let $(\hat{A_1},\hat{\mathcal A_1})$ be an equitably 2-chromatic $P_4$ system of order seven, and $(\hat{A_i},\hat{\mathcal A_i})$ for $2\leq i\leq t$, be an  equitably 2-chromatic $P_4$ system of order six similar to that constructed from Lemma~\ref{lem0.1}. Also, for each $1\leq k<\ell \leq t$, let $\hat{\mathcal A}_{k,\ell}$ be an equitably 2-chromatic decomposition of the edges between $\hat{A_k}$ and $\hat{A_{\ell}}$ into $P_4$ paths which is constructed in a similar manner to Lemma~\ref{lem.h} and Lemma~\ref{lem.i}.

Let $\hat{\mathcal B}=(\bigcup_{i=1}^{t}\hat{\mathcal A_i})\cup \big(\bigcup_{k=1}^{t-1}(\bigcup_{\ell=k+1}^{t}\hat{\mathcal A}_{k,\ell})\big)$, $\hat{R}=\{1,3,\ldots,6t+1\}$, and $\hat{Y}=\{2,4,\ldots,6t\}$. Then $(\hat{V},\hat{\mathcal B})$ is a $P_4(6t+1)$ for which $\hat{R}$ and $\hat{Y}$ constitute an equitable 2-colouring.
\end{Proof}

\begin{Lemma}\label{lem0.4}
Let $t\geq 1$. There exist an equitably 2-chromatic $P_4$ system of order $6t+3$ and a strongly equitably 2-chromatic $P_4$ system of order $6t+4$.
\end{Lemma}
\begin{Proof}
By Lemma~\ref{lem0.3}, there exists a strongly equitably 2-chromatic $P_4(6t)$, $(V,\mathcal B)$, where $V=\{1,\ldots,6t\}$ is the set of points and $\mathcal B$ is the set of blocks. Let $R=\{1,3,\ldots,6t-1\}$ and $Y=\{2,4,\ldots,6t\}$ constitute a strongly equitable 2-colouring of $(V,\mathcal B)$.  Partition the set of points $V$ into $t$ subsets $A_1=\{1,2,3,4,5,6\},\ldots,A_t=\{6t-5,6t-4,6t-3,6t-2,6t-1,6t\}$.

First, we construct an equitably 2-chromatic $P_4(6t+3)$, $(\hat{V},\hat{\mathcal B})$, from $(V,\mathcal B)$. Let $\hat{V}=V\cup\{6t+1,6t+2,6t+3\}$. Decompose the edges between $A_1$ and $\{6t+1,6t+2,6t+3\}$ along with the edges $\{6t+1,6t+2\},\{6t+2,6t+3\},\{6t+3,6t+1\}$ into the set of $P_4$ paths $\hat{\mathcal A_1}$ in a manner similar to Lemma~\ref{lem6}. Also, for each $2\leq i\leq t$, decompose the edges between $A_i$ and $\{6t+1,6t+2,6t+3\}$ into a set $\hat{\mathcal A_i}$ of $P_4$ paths in a manner similar to Lemma~\ref{lem.e}.

Let $\hat{\mathcal B}=\mathcal B\cup (\bigcup_{i=1}^{t}\hat{\mathcal A_i})$, $\hat{R}=\{1,3,\ldots,6t+3\}$, and $\hat{Y}=\{2,4,\ldots,6t+2\}$. Then $(\hat{V},\hat{\mathcal B})$ is a $P_4(6t+3)$ for which $\hat{R}$ and $\hat{Y}$ constitute an equitable 2-colouring.

Next, we construct a strongly equitably 2-chromatic $P_4(6t+4)$, $(\tilde{V},\tilde{\mathcal B})$, from $(V,\mathcal B)$. Let $\tilde{V}=V\cup\{6t+1,6t+2,6t+3,6t+4\}$. Decompose the edges between $A_1$ and $\{6t+1,6t+2,6t+3,6t+4\}$ along with the edges $\{6t+1,6t+2\}$, $\{6t+2,6t+3\}$, $\{6t+3,6t+4\}$, $\{6t+4,6t+1\}$, $\{6t+1,6t+3\}$, $\{6t+2,6t+4\}$ into the set of $P_4$ paths $\tilde{\mathcal A_1}=\big\{(1,6t+1,6t+2,2)$, $(1,6t+3,6t+2,6t+4)$, $(6t+4,2,6t+3,6t+1)$, $(2,6t+1,6t+4,6t+3)$, $(1,6t+2,3,6t+1)$, $(1,6t+4,3,6t+3)$, $(4,6t+1,5,6t+2)$, $(4,6t+3,5,6t+4)$, $(4,6t+2,6,6t+1)$, $(4,6t+4,6,6t+3)\big\}$. Also, for each $2\leq i\leq t$, decompose the edges between $A_i$ and $\{6t+1,6t+2,6t+3,6t+4\}$ into a set $\tilde{\mathcal A_i}$ of $P_4$ paths in a manner similar to Lemma~\ref{lem.f}.

Let $\tilde{\mathcal B}=\mathcal B\cup (\bigcup_{i=1}^{t}\tilde{\mathcal A_i})$, $\tilde{R}=\{1,3,\ldots,6t+3\}$, and $\tilde{Y}=\{2,4,\ldots,6t+4\}$. Then $(\tilde{V},\tilde{\mathcal B})$ is a $P_4(6t+4)$ for which $\tilde{R}$ and $\tilde{Y}$ constitute a strongly equitable 2-colouring.
\end{Proof}

We thus obtain the following theorem.

\begin{Theorem}
For each admissible order $n$, there exists an equitably 2-chromatic $P_4$ system of order $n$.
\end{Theorem}
\begin{Proof}
For $n\leq 7$ apply Lemma~\ref{lem0.1}. Otherwise apply Lemma~\ref{lem0.3} and Lemma~\ref{lem0.4}.
\end{Proof}

\subsection{$\mathbf k$-chromatic $\mathbf {P_4}$ systems for $\mathbf {k\geq 3}$}
In this subsection, we prove that for any integer $k\geq 3$, there exists some integer $n_k$ such that for all admissible $n\geq n_k$, there exists a $k$-chromatic $P_4$ system of order $n$. To prove the main result of this subsection, we first establish several lemmas.

\begin{Lemma}\label{lem9}
Let $t\geq 1$ and $k\geq 3$. If there exists a $k$-chromatic $P_4$ system of order $6t$, then there exist $k$-chromatic $P_4$ systems of order $6t+3$, $6t+4$, $6t+6$, and $6t+7$.
\end{Lemma}
\begin{Proof}
Suppose that there exists a $k$-chromatic $P_4(6t)$, $(V,\mathcal B)$, where $V=\{1,\ldots,6t\}$ is the set of points and $\mathcal B$ is the set of blocks. Let $C_1,C_2,\ldots,C_k$ be the colour classes of a $k$-colouring of $(V,\mathcal B)$. Partition the set of points $V$ into $t$ subsets $A_1=\{1,2,3,4,5,6\},\ldots,A_t=\{6t-5,6t-4,6t-3,6t-2,6t-1,6t\}$.

First, we construct a $k$-chromatic $P_4(6t+3)$, $(V_1,\mathcal B_1)$, from $(V,\mathcal B)$. Let $V_1=V\cup U_1$ where $U_1=\{6t+1,6t+2,6t+3\}$. Let $w_1=6t+1$, $w_2=6t+2$, and $w_3=6t+3$. Decompose the edges between $A_1$ and $U_1$ along with the edges $\{6t+1,6t+2\},\{6t+2,6t+3\},\{6t+3,6t+1\}$ into a set $\mathcal A_1$ of $P_4$ paths in a manner similar to Lemma~\ref{lem6}. For each $2\leq i\leq t$, decompose the edges between $A_i$ and $U_1$ into a set $\mathcal A_i$ of $P_4$ paths  in a manner similar to the 2-chromatic decomposition of $K_{6,3}$ into $P_4$ paths constructed from Lemma~\ref{lem.e}. Let $\mathcal B_1=\mathcal B\cup (\bigcup_{i=1}^{t}\mathcal A_i)$. Then $(V_1,\mathcal B_1)$ is a $P_4(6t+3)$ for which $C_1\cup\{6t+1\}, C_2\cup\{6t+2\},C_3\cup\{6t+3\},C_4,\ldots,C_k$ constitute a $k$-colouring.

Next, we construct a $k$-chromatic $P_4(6t+4)$, $(V_2,\mathcal B_2)$, from $(V,\mathcal B)$. Let $V_2=V\cup U_2$ where $U_2=\{6t+1,6t+2,6t+3,6t+4\}$. Decompose the edges of the complete graph on $U_2$ into a set $\mathcal A$ of $P_4$ paths in a manner similar to the 2-chromatic decomposition of $K_4$ into $P_4$ paths constructed from Lemma~\ref{lem0.1} where $v_1=6t+1$, $v_2=6t+2$, $v_3=6t+3$, and $v_4=6t+4$. For each $1\leq i\leq t$, decompose the edges between $A_i$ and $U_2$ into a set $\mathcal A_i$ of $P_4$ paths in a manner similar to the 2-chromatic decomposition of $K_{6,4}$ into $P_4$ paths constructed from Lemma~\ref{lem.f} where $w_1=6t+1$, $w_2=6t+2$, $w_3=6t+3$, and $w_4=6t+4$. Let $\mathcal B_2=\mathcal B\cup \mathcal A\cup(\bigcup_{i=1}^{t}\mathcal A_i)$. Then $(V_2,\mathcal B_2)$ is a $P_4(6t+4)$ for which $C_1\cup\{6t+1,6t+3\}, C_2\cup\{6t+2,6t+4\},C_3,\ldots,C_k$ constitute a $k$-colouring.

Next, we construct a $k$-chromatic $P_4(6t+6)$, $(V_3,\mathcal B_3)$, from $(V,\mathcal B)$. Let $V_3=V\cup U_3$ where $U_3=\{6t+1,6t+2,6t+3,6t+4,6t+5,6t+6\}$. Let $v_1=6t+1$, $v_2=6t+2$, $v_3=6t+3$, $v_4=6t+4$, $v_5=6t+5$, and $v_6=6t+6$. Decompose the edges of the complete graph on $U_3$ into a set $\mathcal U$ of $P_4$ paths in a manner similar to the 2-chromatic decomposition of $K_6$ into $P_4$ paths constructed from Lemma~\ref{lem0.1}. For each $1\leq i\leq t$, decompose the edges between $A_i$ and $U_3$ into a set $\mathcal A_i$ of $P_4$ paths in a manner similar to the 2-chromatic decomposition of $K_{6,6}$ into $P_4$ paths constructed from Lemma~\ref{lem.h}. Let $\mathcal B_3=\mathcal B\cup \mathcal U\cup(\bigcup_{i=1}^{t}\mathcal A_i)$. Then $(V_3,\mathcal B_3)$ is a $P_4(6t+6)$ for which $C_1\cup\{6t+1,6t+3,6t+5\}, C_2\cup\{6t+2,6t+4,6t+6\},C_3,\ldots,C_k$ constitute a $k$-colouring.

Finally, we construct a $k$-chromatic $P_4(6t+7)$, $(V_4,\mathcal B_4)$, from $(V,\mathcal B)$. Let $V_4=V\cup U_4$ where $U_4=\{6t+1,6t+2,6t+3,6t+4,6t+5,6t+6,6t+7\}$. Decompose the edges of the complete graph on $U_4$ into a set $\mathcal U$ of $P_4$ paths in a manner similar to the 2-chromatic decomposition of $K_7$ into $P_4$ paths constructed from Lemma~\ref{lem0.1} where $v_1=6t+1$, $v_2=6t+2$, $v_3=6t+3$, $v_4=6t+4$, $v_5=6t+5$, $v_6=6t+6$, and $v_7=6t+7$. For each $1\leq i\leq t$, decompose the edges between $A_i$ and $U_4$ into a set $\mathcal A_i$ of $P_4$ paths in a manner similar to the 2-chromatic decomposition of $K_{6,7}$ into $P_4$ paths constructed from Lemma~\ref{lem.i} where $w_1=6t+1$, $w_2=6t+2$, $w_3=6t+3$, $w_4=6t+4$, $w_5=6t+5$, $w_6=6t+6$, and $w_7=6t+7$. Let $\mathcal B_4=\mathcal B\cup \mathcal U\cup(\bigcup_{i=1}^{t}\mathcal A_i)$. Then $(V_4,\mathcal B_4)$ is a $P_4(6t+7)$ for which $C_1\cup\{6t+1,6t+5,6t+7\}, C_2\cup\{6t+2,6t+4,6t+6\},C_3\cup\{6t+3\},\ldots,C_k$ constitute a $k$-colouring.
\end{Proof}

\begin{Lemma}\label{lem10}
Let $t\geq 1$ and $k\geq 3$. If there exists a $k$-chromatic $P_4$ system of order $6t+1$, then there exist $k$-chromatic $P_4$ systems of order $6t+3$, $6t+4$, $6t+6$, and $6t+7$.
\end{Lemma}
\begin{Proof}
Suppose that there exists a $k$-chromatic $P_4(6t+1)$, $(V,\mathcal B)$, where $V=\{1,\ldots,6t+1\}$ is the set of points and $\mathcal B$ is the set of blocks. Let $C_1,C_2,\ldots,C_k$ be the colour classes of a $k$-colouring of $(V,\mathcal B)$. Partition the set of points $V$ into $t$ subsets $A_1=\{1,2,3,4,5,6\}, \ldots, A_{t-1}=\{6t-11,6t-10,6t-9,6t-8,6t-7,6t-6\},A_t=\{6t-5,6t-4,6t-3,6t-2,6t-1,6t,6t+1\}$.

First, we construct a $k$-chromatic $P_4(6t+3)$, $(V_1,\mathcal B_1)$, from $(V,\mathcal B)$. Let $V_1=V\cup U_1$ where $U_1=\{6t+2,6t+3\}$. Decompose the edges between $A_t$ and $U_1$ along with the edge $\{6t+2,6t+3\}$ into a set $\mathcal A_t$ of $P_4$ paths in a manner similar to Lemma~\ref{lem7}. For each $1\leq i\leq t-1$, decompose the edges between $A_i$ and $U_1$ into a set $\mathcal A_i$ of $P_4$ paths in a manner similar to the 2-chromatic decomposition of $K_{6,2}$ into $P_4$ paths constructed from Lemma~\ref{lem.d}. Let $\mathcal B_1=\mathcal B\cup (\bigcup_{i=1}^{t}\mathcal A_i)$. Then $(V_1,\mathcal B_1)$ is a $P_4(6t+3)$ for which $C_1\cup\{6t+2\}, C_2\cup\{6t+3\},C_3,\ldots,C_k$ constitute a $k$-colouring.

Next, we construct a $k$-chromatic $P_4(6t+4)$, $(V_2,\mathcal B_2)$, from $(V,\mathcal B)$. Let $V_2=V\cup U_2$ where $U_2=\{6t+2,6t+3,6t+4\}$. Let $w_1=6t+2$, $w_2=6t+3$, and $w_3=6t+4$. Decompose the edges between $A_t$ and $U_2$ along with the edges $\{6t+2,6t+3\},\{6t+3,6t+4\},\{6t+4,6t+2\}$ into a set $\mathcal A_t$ of $P_4$ paths in a manner similar to Lemma~\ref{lem8}. For each $1\leq i\leq t-1$, decompose the edges between $A_i$ and $U_2$ into a set $\mathcal A_i$ of $P_4$ paths in a manner similar to the 2-chromatic decomposition of $K_{6,3}$ into $P_4$ paths constructed from Lemma~\ref{lem.e}. Let $\mathcal B_2=\mathcal B\cup (\bigcup_{i=1}^{t}\mathcal A_i)$. Then $(V_2,\mathcal B_2)$ is a $P_4(6t+4)$ for which $C_1\cup\{6t+2\}, C_2\cup\{6t+3\},C_3\cup\{6t+4\},C_4\ldots,C_k$ constitute a $k$-colouring.

Next, we construct a $k$-chromatic $P_4(6t+6)$, $(V_3,\mathcal B_3)$, from $(V,\mathcal B)$. Let $V_3=V\cup U_3$ where $U_3=F_1\cup F_2$, $F_1=\{6t+2,6t+3\}$ and $F_2=\{6t+4,6t+5,6t+6\}$. Decompose the edges between $A_t$ and $F_1$ along with the edge $\{6t+2,6t+3\}$ into a set $\mathcal A_t$ of $P_4$ paths in a manner similar to Lemma~\ref{lem7}. For each $1\leq i\leq t-1$, decompose the edges between $A_i$ and $F_1$ into a set $\mathcal A_i$ of $P_4$ paths in a manner similar to the 2-chromatic decomposition of $K_{6,2}$ into $P_4$ paths constructed from Lemma~\ref{lem.d}. Let $w_1=6t+4, w_2=6t+5,$ and $w_3=6t+6$. Decompose the edges between $A_t$ and $F_2$ along with the edges $\{6t+4,6t+5\},\{6t+5,6t+6\},\{6t+6,6t+4\}$ into a set $\mathcal A'_t$ of $P_4$ paths in a manner similar to Lemma~\ref{lem8}. For each $1\leq i\leq t-1$, decompose the edges between $A_i$ and $F_2$ into a set $\mathcal A'_i$ of $P_4$ paths in a manner similar to the 2-chromatic decomposition of $K_{6,3}$ into $P_4$ paths constructed from Lemma~\ref{lem.e}. Decompose the edges between $F_1$ and $F_2$ into a set of $P_4$ paths $\mathcal F=\big\{(6t+4,6t+2,6t+5,6t+3),(6t+4,6t+3,6t+6,6t+2)\big\}$.
Let $\mathcal B_3=\mathcal B\cup (\bigcup_{i=1}^{t}\mathcal A_i)\cup (\bigcup_{i=1}^{t}\mathcal A'_i)\cup \mathcal F$. Then $(V_3,\mathcal B_3)$ is a $P_4(6t+6)$ for which $C_1\cup\{6t+2,6t+5\}, C_2\cup\{6t+3,6t+6\},C_3\cup\{6t+4\},C_4,\ldots,C_k$ constitute a $k$-colouring.

Finally, we construct a $k$-chromatic $P_4(6t+7)$, $(V_4,\mathcal B_4)$, from $(V,\mathcal B)$. Let $V_4=V\cup U_4$ where $U_4=\{6t+2,6t+3,6t+4,6t+5,6t+6,6t+7\}$. Let $v_1=6t+2$, $v_2=6t+3$, $v_3=6t+4$, $v_4=6t+5$, $v_5=6t+6$, and $v_6=6t+7$. Decompose the edges of the complete graph on $U_4$ into a set $\mathcal U$ of $P_4$ paths in a manner similar to the 2-chromatic decomposition of $K_6$ into $P_4$ paths constructed from Lemma~\ref{lem0.1}. For each $1\leq i\leq t-1$, decompose the edges between $A_i$ and $U_4$ into a set $\mathcal A_i$ of $P_4$ paths in a manner similar to the 2-chromatic decomposition of $K_{6,6}$ into $P_4$ paths constructed from Lemma~\ref{lem.h}. Decompose the edges between $A_t$ and $U_4$ into a set $\mathcal A_t$ of $P_4$ paths in a manner similar to the 2-chromatic decomposition of $K_{6,7}$ into $P_4$ paths constructed from Lemma~\ref{lem.i}. Let $\mathcal B_4=\mathcal B\cup \mathcal U\cup(\bigcup_{i=1}^{t}\mathcal A_i)$. Then $(V_4,\mathcal B_4)$ is a $P_4(6t+7)$ for which $C_1\cup\{6t+2,6t+6\}, C_2\cup\{6t+3,6t+5\},C_3\cup\{6t+4,6t+7\},C_4\ldots,C_k$ constitute a $k$-colouring.
\end{Proof}

\begin{Lemma}\label{lem12}
Let $t\geq 1$ and $k\geq 3$. If there exists a $k$-chromatic $P_4$ system of order $6t+3$, then there exist $k$-chromatic $P_4$ systems of order $6t+6$, $6t+7$, $6t+9$, and $6t+10$.
\end{Lemma}
\begin{Proof}
Suppose that there exists a $k$-chromatic $P_4(6t+3)$, $(V,\mathcal B)$, where $V=\{1,\ldots,6t+3\}$ is the set of points and $\mathcal B$ is the set of blocks. Let $C_1,C_2,\ldots,C_k$ be the colour classes of a $k$-colouring of $(V,\mathcal B)$. Partition the set of points $V$ into $t+1$ subsets $A_1=\{1,2,3,4,5,6\}, \ldots,$ $A_t=\{6t-5,6t-4,6t-3,6t-2,6t-1,6t\}$, and $A_{t+1}=\{6t+1,6t+2,6t+3\}$.

First, we construct a $k$-chromatic $P_4(6t+6)$, $(V_1,\mathcal B_1)$, from $(V,\mathcal B)$. Let $V_1=V\cup U_1$ where $U_1=\{6t+4,6t+5,6t+6\}$. Let $w_1=6t+4$, $w_2=6t+5,$ and $w_3=6t+6$. Decompose the edges between $A_1$ and $U_1$ along with the edges of the complete graph on $U_1$ into a set $\mathcal A_1$ of $P_4$ paths in a manner similar to Lemma~\ref{lem6}. For each $2\leq i\leq t$, decompose the edges between $A_i$ and $U_1$ into a set $\mathcal A_i$ of $P_4$ paths in a manner similar to the 2-chromatic decomposition of $K_{6,3}$ into $P_4$ paths constructed from Lemma~\ref{lem.e}. Also, decompose the edges between $A_{t+1}$ and $U_1$ into a set $\mathcal A_{t+1}$ of $P_4$ paths in a manner similar to the 2-chromatic decomposition of $K_{3,3}$ into $P_4$ paths constructed from Lemma~\ref{lem.a}. Let $\mathcal B_1=\mathcal B\cup (\bigcup_{i=1}^{t+1}\mathcal A_i)$. Then $(V_1,\mathcal B_1)$ is a $P_4(6t+6)$ for which $C_1\cup\{6t+4\}, C_2\cup\{6t+5\},C_3\cup\{6t+6\},C_4,\ldots,C_k$ constitute a $k$-colouring.

Next, we construct a $k$-chromatic $P_4(6t+7)$, $(V_2,\mathcal B_2)$, from $(V,\mathcal B)$. Let $V_2=V\cup U_2$ where $U_2=\{6t+4,6t+5,6t+6,6t+7\}$. Decompose the edges of the complete graph on $U_2$ into a set $\mathcal U$ of $P_4$ paths in a manner similar to the equitable 2-chromatic decomposition of $K_4$ into $P_4$ paths constructed from Lemma~\ref{lem0.1} where $v_1=6t+4$, $v_2=6t+5$, $v_3=6t+6$, and $v_4=6t+7$. For each $1\leq i\leq t$, decompose the edges between $A_i$ and $U_2$ into a set $\mathcal A_i$ of $P_4$ paths in a manner similar to the 2-chromatic decomposition of $K_{6,4}$ into $P_4$ paths constructed from Lemma~\ref{lem.f} where $w_1=6t+4$, $w_2=6t+6$, $w_3=6t+5$, and $w_4=6t+7$. Also, decompose the edges between $A_{t+1}$ and $U_2$ into a set $\mathcal A_{t+1}$ of $P_4$ paths in a manner similar to the 2-chromatic decomposition of $K_{4,3}$ into $P_4$ paths constructed from Lemma~\ref{lem.b} where $v_1=6t+4, v_2=6t+5, v_3=6t+6,$ and $v_4=6t+7$. Let $\mathcal B_2=\mathcal B\cup \mathcal U\cup (\bigcup_{i=1}^{t+1}\mathcal A_i)$. Then $(V_2,\mathcal B_2)$ is a $P_4(6t+7)$ for which $C_1\cup\{6t+4,6t+5\}, C_2\cup\{6t+6,6t+7\},C_3,\ldots,C_k$ constitute a $k$-colouring.

Next, we construct a $k$-chromatic $P_4(6t+9)$, $(V_3,\mathcal B_3)$, from $(V,\mathcal B)$. Let $V_3=V\cup U_3$ where $U_3=\{6t+4,6t+5,6t+6,6t+7,6t+8,6t+9\}$. Let $v_1=6t+4$, $v_2=6t+5$, $v_3=6t+6$, $v_4=6t+7$, $v_5=6t+8$, and $v_6=6t+9$. Decompose the edges of the complete graph on $U_3$ into a set $\mathcal U$ of $P_4$ paths in a manner similar to the equitable 2-chromatic decomposition of $K_6$ into $P_4$ paths constructed from Lemma~\ref{lem0.1}. For each $1\leq i\leq t$, decompose the edges between $A_i$ and $U_3$ into a set $\mathcal A_i$ of $P_4$ paths in a manner similar to the 2-chromatic decomposition of $K_{6,6}$ into $P_4$ paths constructed from Lemma~\ref{lem.h}. Also, decompose the edges between $A_{t+1}$ and $U_3$ into a set $\mathcal A_{t+1}$ of $P_4$ paths in a manner similar to the 2-chromatic decomposition of $K_{6,3}$ into $P_4$ paths constructed from Lemma~\ref{lem.e}. Let $\mathcal B_3=\mathcal B\cup \mathcal U\cup (\bigcup_{i=1}^{t+1}\mathcal A_i)$. Then $(V_3,\mathcal B_3)$ is a $P_4(6t+9)$ for which $C_1\cup\{6t+4,6t+6,6t+8\}, C_2\cup\{6t+5,6t+7,6t+9\},C_3,\ldots,C_k$ constitute a $k$-colouring.

Finally, we construct a $k$-chromatic $P_4(6t+10)$, $(V_4,\mathcal B_4)$, from $(V,\mathcal B)$. Let $V_4=V\cup U_4$ where $U_4=\{6t+4,6t+5,6t+6,6t+7,6t+8,6t+9,6t+10\}$. Decompose the edges of the complete graph on $U_4$ into a set $\mathcal U$ of $P_4$ paths in a manner similar to the equitable 2-chromatic decomposition of $K_7$ into $P_4$ paths constructed from Lemma~\ref{lem0.1} where $v_1=6t+4$, $v_2=6t+5$, $v_3=6t+6$, $v_4=6t+7$, $v_5=6t+8$, $v_6=6t+9$, and $v_7=6t+10$. For each $1\leq i\leq t$, decompose the edges between $A_i$ and $U_4$ into a set $\mathcal A_i$ of $P_4$ paths in a manner similar to the 2-chromatic decomposition of $K_{6,7}$ into $P_4$ paths constructed from Lemma~\ref{lem.i} where $w_1=6t+4$, $w_2=6t+5$, $w_3=6t+6$, $w_4=6t+7$, $w_5=6t+8$, $w_6=6t+9$, and $w_7=6t+10$. Also, decompose the edges between $A_{t+1}$ and $U_4$ into a set $\mathcal A_{t+1}$ of $P_4$ paths in a manner similar to the 2-chromatic decomposition of $K_{7,3}$ into $P_4$ paths constructed from Lemma~\ref{lem.c} where $v_1=6t+4$, $v_2=6t+5$, $v_3=6t+6$, $v_4=6t+7$, $v_5=6t+8$, $v_6=6t+9$, and $v_7=6t+10$. Let $\mathcal B_4=\mathcal B\cup \mathcal U\cup (\bigcup_{i=1}^{t+1}\mathcal A_i)$. Then $(V_4,\mathcal B_4)$ is a $P_4(6t+10)$ for which $C_1\cup\{6t+4,6t+8\}, C_2\cup\{6t+5,6t+7,6t+9\},C_3\cup\{6t+6,6t+10\},C_4,\ldots,C_k$ constitute a $k$-colouring.
\end{Proof}

\begin{Lemma}\label{lem16}
Let $t\geq 1$ and $k\geq 3$. If there exists a $k$-chromatic $P_4$ system of order $6t+4$, then there exist $k$-chromatic $P_4$ systems of order $6t+6$, $6t+7$, $6t+9$, and $6t+10$.
\end{Lemma}
\begin{Proof}
Suppose that there exists a $k$-chromatic $P_4(6t+4)$, $(V,\mathcal B)$, where $V=\{1,\ldots,6t+4\}$ is the set of points and $\mathcal B$ is the set of blocks. Let $C_1,C_2,\ldots,C_k$ be the colour classes of a $k$-colouring of $(V,\mathcal B)$. Partition the set of points $V$ into $t+1$ subsets $A_1=\{1,2,3,4,5,6\}, \ldots, A_t=\{6t-5,6t-4,6t-3,6t-2,6t-1,6t\}$, and $A_{t+1}=\{6t+1,6t+2,6t+3,6t+4\}$.

First, we construct a $k$-chromatic $P_4(6t+6)$, $(V_1,\mathcal B_1)$, from $(V,\mathcal B)$. Let $V_1=V\cup U_1$ where $U_1=\{6t+5,6t+6\}$.  For each $1\leq i\leq t$, decompose the edges between $A_i$ and $U_1$ into a set  $\mathcal A_i$ of $P_4$ paths in a manner similar to the 2-chromatic decomposition of $K_{6,2}$ into $P_4$ paths constructed from Lemma~\ref{lem.d}. Decompose the edges between $A_{t+1}$ and $U_1$ along with the edges of the complete graph on $U_1$ into a set $\mathcal A_{t+1}$ of $P_4$ paths in a manner similar to Lemma~\ref{lem13}.  Let $\mathcal B_1=\mathcal B\cup (\bigcup_{i=1}^{t+1}\mathcal A_i)$. Then $(V_1,\mathcal B_1)$ is a $P_4(6t+6)$ for which $C_1\cup\{6t+5\}, C_2\cup\{6t+6\},C_3,\ldots,C_k$ constitute a $k$-colouring.

Next, we construct a $k$-chromatic $P_4(6t+7)$, $(V_2,\mathcal B_2)$, from $(V,\mathcal B)$. Let $V_2=V\cup U_2$ where $U_2=\{6t+5,6t+6,6t+7\}$.  Let $w_1=6t+5$, $w_2=6t+6$, and $w_3=6t+7$. Decompose the edges between $A_1$ and $U_2$ along with the edges of the complete graph on $U_2$ into a set $\mathcal A_1$ of $P_4$ paths in a manner similar to Lemma~\ref{lem6}. For each $2\leq i\leq t$, decompose the edges between $A_i$ and $U_2$ into a set $\mathcal A_i$ of $P_4$ paths in a manner similar to the 2-chromatic decomposition of $K_{6,3}$ into $P_4$ paths constructed from Lemma~\ref{lem.e}. Also, decompose the edges between $A_{t+1}$ and $U_2$ into a set $\mathcal A_{t+1}$ of $P_4$ paths in a manner similar to the 2-chromatic decomposition of $K_{4,3}$ into $P_4$ paths constructed from Lemma~\ref{lem.b}. Let $\mathcal B_2=\mathcal B\cup (\bigcup_{i=1}^{t+1}\mathcal A_i)$. Then $(V_2,\mathcal B_2)$ is a $P_4(6t+7)$ for which $C_1\cup\{6t+5\}, C_2\cup\{6t+6\},C_3\cup\{6t+7\},C_4,\ldots,C_k$ constitute a $k$-colouring.

Next, we construct a $k$-chromatic $P_4(6t+9)$, $(V_3,\mathcal B_3)$, from $(V,\mathcal B)$. Let $V_3=V\cup U_3$ where $U_3=\{6t+5,6t+6,6t+7,6t+8,6t+9\}$.  Let $w_1=6t+5$, $w_2=6t+6$, $w_3=6t+7$, $w_4=6t+8$, and $w_5=6t+9$. For each $1\leq i\leq t$, decompose the edges between $A_i$ and $U_3$ into a set $\mathcal A_i$ of $P_4$ paths in a manner similar to the 2-chromatic decomposition of $K_{6,5}$ into $P_4$ paths constructed from Lemma~\ref{lem.g}. Decompose the edges between $A_{t+1}$ and $U_3$ along with the edges of the complete graph on $U_3$ into a set $\mathcal A_{t+1}$ of $P_4$ paths in a manner similar to Lemma~\ref{lem14}.  Let $\mathcal B_3=\mathcal B\cup (\bigcup_{i=1}^{t+1}\mathcal A_i)$. Then $(V_3,\mathcal B_3)$ is a $P_4(6t+9)$ for which $C_1\cup\{6t+5,6t+8\}, C_2\cup\{6t+6,6t+9\},C_3\cup\{6t+7\},C_4,\ldots,C_k$ constitute a $k$-colouring.

Finally, we construct a $k$-chromatic $P_4(6t+10)$, $(V_4,\mathcal B_4)$, from $(V,\mathcal B)$. Let $V_4=V\cup U_4$ where $U_4=\{6t+5,6t+6,6t+7,6t+8,6t+9,6t+10\}$. Let $v_1=6t+5$, $v_2=6t+6$, $v_3=6t+7$, $v_4=6t+8$, $v_5=6t+9$, and $v_6=6t+10$. Decompose the edges of the complete graph on $U_4$ into a set $\mathcal U$ of $P_4$ paths in a manner similar to the equitable 2-chromatic decomposition of $K_6$ into $P_4$ paths constructed from Lemma~\ref{lem0.1}. For each $1\leq i\leq t$, decompose the edges between $A_i$ and $U_4$ into a set $\mathcal A_i$ of $P_4$ paths in a manner similar to the 2-chromatic decomposition of $K_{6,6}$ into $P_4$ paths constructed from Lemma~\ref{lem.h}. Also, decompose the edges between $A_{t+1}$ and $U_4$ into a set $\mathcal A_{t+1}$ of $P_4$ paths in a manner similar to the 2-chromatic decomposition of $K_{4,6}$ into $P_4$ paths constructed from Lemma~\ref{lem.f}. Let $\mathcal B_4=\mathcal B\cup \mathcal U\cup (\bigcup_{i=1}^{t+1}\mathcal A_i)$. Then $(V_4,\mathcal B_4)$ is a $P_4(6t+10)$ for which $C_1\cup\{6t+5,6t+9\}, C_2\cup\{6t+6,6t+8\},C_3\cup\{6t+7,6t+10\},C_4,\ldots,C_k$ constitute a $k$-colouring.
\end{Proof}

We thus obtain the following theorem.

\begin{Theorem}
For any integer $k\geq 3$, there exists some integer $n_k$ such that for all admissible $n\geq n_k$, there exists a $k$-chromatic $P_4$ system of order $n$.
\end{Theorem}
\begin{Proof}
By Corollary~\ref{cor1}, for any integer $k\geq 3$, there exists a $k$-chromatic $P_4$ system of some order $n_k$. Apply
Lemma~\ref{lem9}, Lemma~\ref{lem10}, Lemma~\ref{lem12}, and Lemma~\ref{lem16} to construct a $k$-chromatic $P_4$ system of order $n$ for all admissible $n\geq n_k$.
\end{Proof}

\section{Unique colourings of $P_4$ systems}
We now investigate uniquely 2-chromatic $P_4$ systems. We commence by showing that there exists a strongly equitably 2-chromatic  $P_4$ system for which two particular vertices are forced to have different colours. Note that we denote the colour of vertex $v\in V$ by $c(v)$.

\begin{Lemma}\label{lem4.1}
There exists a strongly equitably 2-chromatic $P_4$ system $(V,\mathcal B)$ of order 28 with two specific vertices that cannot share a colour in any 2-colouring.
\end{Lemma}

\begin{Proof}
We construct a strongly equitably 2-chromatic $P_4$ system $(V,\mathcal B)$ of order 28 with colour classes $C_w$ and $C_b$ such that each vertex of $C_w$ has colour $w$ and each vertex of $C_b$ has colour $b$. Let $V=\{1,2,3,\ldots,27,28\}$. We will specify a block set $\mathcal B$ such that $c(27)$ cannot equal $c(28)$. By way of contradiction, suppose that $c(27)=c(28)=w$. Let $\mathcal B_1=\bigcup\limits_{i=1}^{12}\{(27,2i-1,2i,28)\}$. Then for each $i\in \{1,\ldots,12\}$, $b\in \{c(2i-1),c(2i)\}$ for otherwise the block $(27,2i-1,2i,28)$ would be monochromatic.

For distinct values $i,j,\ell \in \{1,2,\ldots,12\}$, let $\mathcal B(i,j,\ell)=\{(2i-1,2j,2\ell,2i),(2i-1,2j-1,2\ell-1,2i),(2i-1,2\ell -1,2j,2i),(2i-1,2\ell,2j-1,2i)\}$. If $c(2i-1)=c(2i)=b$, then one of the paths of $\mathcal B(i,j,\ell)$ is monochromatic as $b\in \{c(2j-1), c(2j)\}$ and $b\in \{c(2\ell-1), c(2\ell)\}$. Hence $c(2i-1)\neq c(2i)$ and $\{c(2i-1),c(2i)\}=\{b,w\}$.
Let $\mathcal B_2=\big\{\mathcal B(1,2,3)$, $\mathcal B(4,5,6)$, $\mathcal B(7,8,9)$, $\mathcal B(10,11,12)\big\}$
and observe that $c(1)\neq c(2)$, $c(7)\neq c(8)$, $c(13)\neq c(14)$, and $c(19)\neq c(20)$.

Let $\mathcal B_3=\big\{(7,1,13,20)$, $(1,19,14,7)$, $(7,2,14,20)$, $(2,20,7,13)$, $(13,8,1,20)$, $(2,8,19,13)$, $(1,14,8,20)$, $(7,19,2,13)\big\}$. Then $c(1)=c(7)=c(13)=c(19)$ and $c(2)=c(8)=c(14)=c(20)$ for otherwise there would be a monochromatic path in $\mathcal B_3$.

Let $\mathcal B_4=\big\{\mathcal B(2,4,9)$, $\mathcal B(5,1,12)$, $\mathcal B(8,3,10)$, $\mathcal B(11,6,7)\big\}$ and observe that $c(3)\neq c(4)$, $c(9)\neq c(10)$, $c(15)\neq c(16)$, and $c(21)\neq c(22)$.

Let $\mathcal B_{5}=\big\{(9,3,15,22)$, $(3,21,16,9)$, $(9,4,16,22)$, $(4,22,9,15)$, $(15,10,3,22)$, $(4,10,21,15)$, $(3,16,10,22)$, $(9,21,4,15)\big\}$. Then $c(3)=c(9)=c(15)=c(21)$ and $c(4)=c(10)=c(16)=c(22)$ for otherwise there would be a monochromatic path in $\mathcal B_{5}$.

Let $\mathcal B_{c,d}(a_1,a_2,b_1,b_2)=\big\{(c,a_1,b_1,a_2), (d,a_2,b_2,a_1)\big\}$. Then let
$\mathcal B_{6}=\mathcal B_{25,26}(1,2,11,12)$ $\cup$ $\mathcal B_{28,27}(1,2,15,16)$ $\cup$ $\mathcal B_{26,25}(1,2,17,18)$ $\cup$ $\mathcal B_{28,27}(21,22,1,2)$ $\cup$
$\mathcal B_{25,26}(3,4,11,12)$ $\cup$
$\mathcal B_{26,25}(3,4,23,24)$ $\cup$
$\mathcal B_{28,27}(5,6,7,8)$ $\cup$ $\mathcal B_{25,26}(5,6,9,10)$ $\cup$ $\mathcal B_{26,25}(5,6,11,12)$ $\cup$ $\mathcal B_{26,25}(13,14,5,6)$ $\cup$ $\mathcal B_{28,27}(17,18,5,6)$ $\cup$ $\mathcal B_{25,26}(21,22,5,6)$ $\cup$ $\mathcal B_{26,25}(23,24,5,6)$ $\cup$
$\mathcal B_{28,27}(7,8,15,16)$ $\cup$ $\mathcal B_{25,26}(7,8,21,22)$ $\cup$ $\mathcal B_{26,25}(7,8,23,24)$ $\cup$
$\mathcal B_{26,25}(9,10,13,14)$ $\cup$ $\mathcal B_{25,26}(9,10,17,18)$ $\cup$
$\mathcal B_{28,27}(11,12,15,16)$ $\cup$ $\mathcal B_{25,26}(11,12,17,18)$ $\cup$ $\mathcal B_{26,25}(11,12,19,20)$ $\cup$ $\mathcal B_{25,26}(23,24,11,12)$ $\cup$
$\mathcal B_{25,26}(13,14,23,24)$ $\cup$ $\mathcal B_{28,27}(15,16,23,24)$ $\cup$
$\mathcal B_{25,26}(17,18,19,20)$ $\cup$ $\mathcal B_{26,25}(21,22,17,18)$ $\cup$ $\mathcal B_{26,25}(17,18,23,24)$.

Let
$\mathcal B_{7}=\big\{(15,25,26,16)$, $(15,26,19,25)$, $(16,25,20,26)$, $(28,26,27,14)$, $(13,28,25,27)$, $(23,28,27,24)$,
$(14,3,28,19)$, $(13,4,19,10)$, $(14,4,20,27)$, $(3,19,9,28)$, $(4,27,10,20)$, $(13,3,20,9)\big\}$.

Let ${\mathcal B}=\bigcup\limits_{i=1}^{7}\mathcal B_{i}$ and observe that $(V,{\mathcal B})$ is a $P_{4}(28)$.



Recall that $c(13)=c(19)$ and $c(4)=c(10)$. Observe that $c(13)\neq c(4)$ for otherwise the block $(13,4,19,10)$ from $\mathcal B_{7}$ would be monochromatic. If $c(13)=w$, then $c(3)=c(9)=w$ which yields to the block $(3,19,9,28)$ from $\mathcal B_{7}$ being monochromatic. If $c(13)=b$, then $c(14)=c(20)=w$ which yields to the block $(14,4,20,27)$ from $\mathcal B_{7}$ being monochromatic. We have thus established a contradiction and hence for any 2-colouring of $(V,\mathcal B)$ it must be that $c(27)\neq c(28)$. Finally observe that $C_w=\{1,3,5,\ldots,27\}$ and $C_b=\{2,4,6,\ldots,28\}$ exhibit a strongly equitable 2-colouring of $(V,\mathcal B)$.
\end{Proof}


We now show that there exists a uniquely 2-chromatic $P_4$ system of order 109.

\begin{Theorem}\label{Thm-109}
There exists a uniquely 2-chromatic $P_4$ system of order 109.
\end{Theorem}
\begin{Proof}
We construct a uniquely 2-chromatic $P_4$ system $(P,\mathcal P)$ of order 109 with colour classes $C_w$ and $C_b$ such that each vertex of $C_w$ has colour $w$ and each vertex of $C_b$ has colour $b$. Let $P=\{1,\ldots,28\}$ $\cup$ $\{1',\ldots,27'\}$ $\cup$ $\{1'',\ldots,27''\}$ $\cup$ $\{1''',\ldots,27'''\}$.

Let $(V,\mathcal B)$ be the strongly equitably 2-chromatic $P_4(28)$ constructed in the proof of Lemma~\ref{lem4.1}. Interchange the names of points 27 and 28 to yield an isomorphic $P_4(28)$ with point set $V$ and block set $\hat{\mathcal B}$ such that points 27 and 28 cannot have the same colour in any 2-colouring of $(V,\hat{\mathcal B})$ and the partition of $V$ given by the sets $\{1,3,5,\ldots,25\} \cup \{28\}$ and $\{2,4,6,\ldots,26\}\cup \{27\}$ is a valid 2-colouring. Now let $(V',\mathcal B')$ be obtained from $(V,\hat{\mathcal B})$ by placing a single prime on each of the points of $\{1,\ldots,27\}$. So $V'=\{1',2',3',\ldots,27',28\}$ and blocks such as $(4,28,10,20)$ of $\hat{\mathcal B}$ yield blocks such as $(4',28,10',20')$. Also form $(V'',\mathcal B'')$ and $(V''',\mathcal B''')$ in a similar manner. Now observe that 28 cannot share a colour with any of $27, 27', 27''$ or $27'''$ in any 2-colouring of any design that contains the blocks of $\hat{\mathcal B} \cup \mathcal B' \cup \mathcal B'' \cup \mathcal B'''$. Without loss of generality, let $c(28)=w$. Then $c(27)=c(27')=c(27'')=c(27''')=b$.

Let $\mathcal P_1=\{(27,1',27'',27')$, $(27,3',27'',27''')$, $(27,1'',27''',27')$, $(27',27,1''',27'')\}$. Then $c(1')=c(3')=c(1'')=c(1''')=w$. For each $1\leq k\leq 12$, let $\mathcal C_k=\{((2k)'',1',(2k-1)'',3')$, $(27,(2k+1)'',27',(2k)'')\}$ and $\mathcal D_k=\{((2k)''',1',(2k-1)''',3')$, $(27,(2k+1)''',27',(2k)''')\}$. Also, let $\mathcal P_2=\{(26'',1',25'',3')$, $(26''',1',25''',3')\}$. Since $c(1')=c(3')=c(1'')=w$, then for $k=1$ it follows that $c(3'')=w$ and $c(2'')=b$. Similarly, by iteration for each $1\leq k\leq 12$, $c((2k+1)'')=w$ and $c((2k)'')=b$. Also $c(26'')=b$. Likewise, for each $1\leq k\leq 13$, $c((2k-1)''')=w$ and $c((2k)''')=b$.

Let $\mathcal F_k=\{(26'',(2k+1)',(2k)''',27'')\}$, $2\leq k\leq 12$. Then for each $2\leq k\leq 12$, $c((2k+1)')=w$. Also, let $\mathcal H_k=\{(1''',(2k)',3''',(2k-1)''\}$, $1\leq k\leq 13$. Then for each $1\leq k\leq 13$, $c((2k)')=b$.

Let $\mathcal M_k=\{(2''',2k-1,(2k)'',27''')\}$, $1\leq k\leq 13$. Then for each $1\leq k\leq 13$, $c(2k-1)=w$. Also, let $\mathcal N_k=\{(2k-1,1',2k,1'')\}$, $1\leq k\leq 13$. Then for each $1\leq k\leq 13$, $c(2k)=b$.

At this point every vertex of $P$ is coloured, and the colouring is unique (up to a permutation of colours).
We now decompose the remaining edges into $P_4$ paths such that each pair of distinct elements of $P$ belongs to at most one $P_4$ path.
These remaining $P_4$ paths are non-critical in the sense that the uniqueness of the colouring is forced by the blocks that have
already been included in the system.

For each $k\in \{1,3,5,\ldots,25\}$,
let $\mathcal R_k=\{(27',k,2',k+1)$, $(26',k,4',k+1)$, $(25',k,6',k+1)$, $(24',k,8',k+1)$, $(23',k,10',k+1)$, $(22',k,12',k+1)$, $(21',k,14',k+1)$, $(20',k,16',k+1)$, $(19',k,18',k+1)$, $(27',k+1,3',k)$, $(26',k+1,5',k)$, $(25',k+1,7',k)$, $(24',k+1,9',k)$, $(23',k+1,11',k)$, $(22',k+1,13',k)$, $(21',k+1,15',k)$, $(20',k+1,17',k)\}$.

For each $k \in \{1,3,5,\ldots,25\}$,
let $\mathcal S_k =  \{(27'',k,2'',k+1)$, $(26'',k,4'',k+1)$, $(25'',k,6'',k+1)$, $(24'',k,8'',k+1)$, $(23'',k,10'',k+1)$, $(22'',k,12'',k+1)$, $(21'',k,14'',k+1)$, $(20'',k,16'',k+1)$, $(1'',k,19'',k+1)\}$.
Let $\mathcal S = \{(18'',25,4'',26)$, $(18'',23,8'',24)$, $(18'',21,12'',22)$, $(18'',19,16'',20)\}$
$\cup$ $\big( \bigcup_{k \in \{1,3,5,\ldots,25\}} \mathcal S_k$ $\setminus$
$\{(27'',1,2'',2)$, $(26'',3,4'',4)$, $(25'',5,6'',6)$, $(24'',7,8'',8)$, $(23'',9,10'',10)$, $(22'',11,12'',12)$, $(21'',13,14'',14)$, $(20'',15,16'',16)$, $(26'',25,4'',26)$, $(24'',23,8'',24)$, $(22'',21,12'',22)$, $(20'',19,16'',20)\} \big)$.

For each $k \in \{1,3,5,\ldots,25\}$,
let
$\mathcal T_k^1$ $=\{(27'',k+1,3'',k)$, $(26'',k+1,5'',k)$, $(25'',k+1,7'',k)$, $(24'',k+1,9'',k)$, $(23'',k+1,11'',k)$, $(22'',k+1,13'',k)$, $(21'',k+1,15'',k)$, $(20'',k+1,17'',k)\}$,
let $\mathcal T_k^2=\{(27''',k,1''',k+1)$, $(26''',k,4''',k+1)$, $(25''',k,6''',k+1)$, $(24''',k,8''',k+1)$, $(23''',k,10''',k+1)$, $(22''',k,12''',k+1)$, $(21''',k,14''',k+1)$, $(20''',k,16''',k+1)\}$ and
$\mathcal T_k^3=\{(27''',k+1,3''',k)$, $(26''',k+1,5''',k)$, $(25''',k+1,7''',k)$, $(24''',k+1,9''',k)$, $(23''',k+1,11''',k)$, $(22''',k+1,13''',k)$, $(21''',k+1,15''',k)$, $(20''',k+1,17''',k)$, $(2''',k+1,18''',k)\}$.

For each $k\in \{1,3,5,\ldots,23\}$,
let $\mathcal U_k=\{(26',k'',2',(k+1)'')$, $(25',k'',5',(k+1)'')$, $(24',k'',7',(k+1)'')$, $(23',k'',9',(k+1)'')$, $(22',k'',11',(k+1)'')$, $(21',k'',13',(k+1)'')$, $(20',k'',15',(k+1)'')$, $(19',k'',17',(k+1)'')$,
$(26',(k+1)'',4',k'')$, $(25',(k+1)'',6',k'')$, $(24',(k+1)'',8',k'')$, $(23',(k+1)'',10',k'')$, $(22',(k+1)'',12',k'')$, $(21',(k+1)'',14',k'')$, $(20',(k+1)'',16',k'')$, $(19',(k+1)'',18',k'')\}$.
Let $\mathcal U =
\{(26',25'',2',26'')$, $(25',25'',5',27'')$, $(24',25'',7',27'')$, $(23',25'',9',27'')$, $(22',25'',11',27'')$, $(21',25'',13',27'')$, $(20',25'',15',27'')$, $(19',25'',17',27'')$,
$(26',26'',4',25'')$,
$(24',26'',8',25'')$,
$(22',26'',12',25'')$,
$(20',26'',16',25'')$,
$(19',27'',18',25'')$, $(21',27'',14',25'')$, $(23',27'',10',25'')$, $(25',27'',6',25'')\}$
$\cup$
$(\bigcup_{k \in \{1,3,5,\ldots,23\}} \mathcal U_k)$.

For each $k\in \{5,7,9,\ldots,25\}$, let
$\mathcal V_k=\{(26',k''',2',(k+1)''')$, $(25',k''',5',(k+1)''')$, $(24',k''',7',(k+1)''')$, $(23',k''',9',(k+1)''')$, $(22',k''',11',(k+1)''')$, $(21',k''',13',(k+1)''')$, $(20',k''',15',(k+1)''')$, $(19',k''',17',(k+1)''')\}$.
Let $\mathcal V =
\{(24',5''',7',27''')$, $(23',7''',9',27''')$, $(22',9''',11',27''')$, $(21',11''',13',27''')$, $(20',13''',15',27''')$, $(19',15''',17',27''')\}$
$\cup$ $\big( \bigcup_{k \in \{5,7,9,\ldots,25\}} \mathcal V_k$ $\setminus$
$\{(24',5''',7',6''')$, $(23',7''',9',8''')$, $(22',9''',11',10''')$, $(21',11''',13',12''')$, $(20',13''',15',14''')$, $(19',15''',17',16''')\}\big)$.

For each $k\in \{5,7,9,\ldots,25\}$, let
$\mathcal W_k=\{(26',(k+1)''',4',k''')$, $(25',(k+1)''',6',k''')$, $(24',(k+1)''',8',k''')$, $(23',(k+1)''',10',k''')$, $(22',(k+1)''',12',k''')$, $(21',(k+1)''',14',k''')$, $(20',(k+1)''',16',k''')$, $(19',(k+1)''',18',k''')\}$.
Let $\mathcal W =
\{(23''',6',27''',25')$, $(21''',10',27''',23')$, $(19''',14',27''',21')$, $(17''',18',27''',19')\}$
$\cup$ $\big( \bigcup_{k \in \{5,7,9,\ldots,25\}} \mathcal W_k$ $\setminus$
$\{(25',24''',6',23''')$, $(23',22''',10',21''')$, $(21',20''',14',19''')$, $(19',18''',18',17''')\}\big)$.

For each $k\in \{1,3,5,\ldots,25\}$, let $\mathcal X_k=\{(26''',k'',1''',(k+1)'')$, $(25''',k'',4''',(k+1)'')$, $(24''',k'',6''',(k+1)'')$, $(23''',k'',8''',(k+1)'')$, $(22''',k'',10''',(k+1)'')$, $(21''',k'',12''',(k+1)'')$, $(20''',k'',14''',(k+1)'')$, $(19''',k'',16''',(k+1)'')$,
$(26''',(k+1)'',2''',k'')$, $(25''',(k+1)'',5''',k'')$, $(24''',(k+1)'',7''',k'')$, $(23''',(k+1)'',9''',k'')$, $(22''',(k+1)'',11''',k'')$,
$(21''',(k+1)'',13''',k'')$, $(20''',(k+1)'',15''',k'')$, $(19''',(k+1)'',17''',k'')$, $(3''',(k+1)'',18''',k'')\}$.

Let $\mathcal Y$ be the set of 94 $P_4$ paths listed in Table~\ref{Table-94blocks}.
This collection of paths was found by performing a computer search for a set of paths
that decomposes those edges which had not yet been used,
such that each path has at least one vertex of each colour.

Let $\mathcal P=\hat{\mathcal B} \cup \mathcal B' \cup \mathcal B'' \cup \mathcal B'''$ $\cup$ $(\bigcup \limits_{k=1}^{2}\mathcal P_k)$ $\cup$ $(\bigcup \limits_{k=1}^{12}\mathcal C_k)$ $\cup$ $(\bigcup \limits_{k=1}^{12}\mathcal D_k)$ $\cup$ $(\bigcup \limits_{k=2}^{12}\mathcal F_k)$ $\cup$ $(\bigcup \limits_{k=1}^{13}\mathcal H_k)$ $\cup$ $(\bigcup \limits_{k=1}^{13}\mathcal M_k)$ $\cup$ $(\bigcup \limits_{k=1}^{13}\mathcal N_k)$ $\cup$ $\big(\bigcup\limits_{k\in\{1,3,5,\ldots,25\}}\mathcal R_k\big)$
$\cup$ $\mathcal S$
$\cup$ $\big(\bigcup\limits_{k\in\{1,3,5,\ldots,25\}}(\mathcal T_k^1 \cup \mathcal T_k^2 \cup \mathcal T_k^3)\big)$
$\cup$ $\mathcal U$
$\cup$ $\mathcal V$
$\cup$ $\mathcal W$
$\cup$ $\big(\bigcup\limits_{k\in\{1,3,5,\ldots,25\}}\mathcal X_k \big)$
$\cup$ $\mathcal Y$.
Then $(P,\mathcal P)$ is a uniquely 2-chromatic $P_4(109)$ for which $C_w=\{1,3,5,\ldots,25\}$ $\cup$ $\{1',3',5',\ldots,25'\}$ $\cup$ $\{1'',3'',5'',\ldots,25''\}$ $\cup$ $\{1''',3''',5''',\ldots,25'''\}$ $\cup$ $\{28\}$ and $C_b=\{2,4,6,\ldots,26\}$ $\cup$ $\{2',4',6',\ldots,26'\}$ $\cup$ $\{2'',4'',6'',\ldots,26''\}$ $\cup$ $\{2''',4''',6''',\ldots,26'''\}$ $\cup$ $\{27,27',27'',27'''\}$ constitute a unique 2-colouring.
\end{Proof}

\begin{table}[htbp]
\begin{center}
\begin{tabular}{cccc}
$(17,19''',1,18'')$&
$(1,27'',27,2')$&
$(19,19''',2,19')$&
$(21,19''',3,18'')$
\\
$(3,26'',27,4')$&
$(23,19''',4,19')$&
$(25,19''',5,18'')$&
$(5,25'',27''',27)$
\\
$(1',27''',2',27'')$&
$(1'',27',26'',3')$&
$(3'',27''',3',2'')$&
$(5'',27''',4',27'')$
\\
$(7'',27''',5',27)$&
$(9'',27''',8',27)$&
$(11'',27''',12',27)$&
$(13'',27''',16',27)$
\\
$(15'',27''',20',27)$&
$(17'',27''',22',27)$&
$(19'',27''',24',27)$&
$(5''',27'',8',2''')$
\\
$(8',4''',27,7')$&
$(7''',27'',12',2''')$&
$(12',4''',3',4'')$&
$(9''',27'',16',2''')$
\\
$(16',4''',7',1''')$&
$(11''',27'',20',2''')$&
$(20',4''',9',27)$&
$(13''',27'',22',2''')$
\\
$(22',4''',11',27)$&
$(15''',27'',24',2''')$&
$(24',4''',13',27)$&
$(17''',27'',26',27)$
\\
$(21''',27'',2''',27)$&
$(23''',27'',3''',5')$&
$(25''',27'',19''',6)$&
$(27'',26''',27,15')$
\\
$(2,2'',27,17')$&
$(2,18'',6,19')$&
$(6,6'',27,19')$&
$(6'',3',8'',8)$
\\
$(8'',27,21',1''')$&
$(4,4'',27,23')$&
$(4,18'',7,24'')$&
$(7,19''',8,19')$
\\
$(8,18'',9,23'')$&
$(9,19''',10,19')$&
$(23'',27''',26',2''')$&
$(26',4''',15',1''')$
\\
$(27''',21'',13,18'')$&
$(13,19''',11,18'')$&
$(11,22'',27,6')$&
$(22'',3',10'',10)$
\\
$(10,18'',12,19')$&
$(10'',27,25',1''')$&
$(12,12'',3',14'')$&
$(12,19''',14,19')$
\\
$(12'',27,16'',3')$&
$(16'',16,19',18)$&
$(14'',14,18'',15)$&
$(14'',27,18'',3')$
\\
$(15,20'',27,10')$&
$(15,19''',16,18'')$&
$(20'',3',24'',27)$&
$(18,18'',20,19')$
\\
$(18,19''',22,19')$&
$(20,19''',24,18'')$&
$(22,18'',26,19''')$&
$(24,19',1''',5')$
\\
$(26,19',2''',2')$&
$(2',4''',17',1''')$&
$(5',2''',3',6''')$&
$(6''',27,8''',3')$
\\
$(4',2''',7',3''')$&
$(4',4''',19',3''')$&
$(15',2''',6',26'')$&
$(15',3''',9',2''')$
\\
$(9',1''',11',2''')$&
$(11',3''',13',2''')$&
$(13',1''',23',2''')$&
$(1''',27',26''',3')$
\\
$(4''',6',24''',3')$&
$(24''',27,10''',3')$&
$(17',2''',10',26'')$&
$(17',3''',21',2''')$
\\
$(21',4''',10',22''')$&
$(23',3''',25',2''')$&
$(23',4''',14',26'')$&
$(25',4''',18',26'')$
\\
$(2''',14',20''',3')$&
$(14',27,12''',3')$&
$(2''',18',18''',3')$&
$(18',27,14''',3')$
\\
$(18''',27,16''',3')$&
$(3',22''',27,20''')$
\end{tabular}
\end{center}
\caption{The set $\mathcal Y$ of $P_4$ paths used in the proof of Theorem~\ref{Thm-109}.}
\label{Table-94blocks}
\end{table}

We now present several results that collectively establish the existence of uniquely 2-chromatic $P_4$ systems of many orders.

\begin{Lemma}\label{cong1_2}
Let $(P,\mathcal P)$ be a uniquely 2-chromatic $P_4$ system of order $n$ where $n \equiv 0,1,3,4$ (mod 6), with colour classes $C_w$ and $C_b$ such that $|C_w| \geq 6$ and $|C_b| \geq 6$. If there exist  non-critical $P_4$ paths $(a_6,a_2,a_3,a_5)$ and $(a_2,a_4,a_1,a_3)$ in $\mathcal P$ such that $c(a_1)=c(a_3)=c(a_5)=w$ and $c(a_2)=c(a_4)=c(a_6)=b$, then there exists a uniquely 2-chromatic $P_4$ system of order $n+6$.
\end{Lemma}

\begin{Proof}
Let $\{a_1,a_3,a_5,a_7,a_9,a_{11}\} \subseteq C_w$ and $\{a_2,a_4,a_6,a_8,a_{10},a_{12}\} \subseteq C_b$. Also, let $P'=P\cup \{v_1,v_2,v_3,v_4,v_5,v_6\}$. We will remove the two non-critical paths $(a_6,a_2,a_3,a_5)$ and $(a_2,a_4,a_1,a_3)$ from $\mathcal P$, so that we can reuse their edges in a different manner.

Let $\mathcal A= \{(a_1,v_4,a_3,a_5)$, $(v_4,a_4,v_5,a_8)$, $(v_4,a_{10},v_3,a_{12})$,
$(a_4, a_2,v_1, a_6)$, $(a_5,v_1,a_7,v_6)$, $(v_2,a_9,v_1,a_{11})\}$. Then $c(v_1)=c(v_3)=c(v_5)=w$ and $c(v_2)=c(v_4)=c(v_6)=b$.  Decompose the remaining edges between $\{a_1,a_3\}$ and $\{v_1,v_2,v_3,v_4,v_5,v_6\}$ along with the edges $(a_1,a_3)$ and $(a_2,a_3)$ into the set of $P_4$ paths
$\mathcal B = \{(a_2,a_3,a_1,v_5)$, $(a_1,v_6,a_3,v_5)$, $(a_1,v_1,a_3,v_2)$, $(a_3,v_3,a_1,v_2)\}$.
Decompose the remaining edges between $\{a_5,a_7\}$ and $\{v_1,v_2,v_3,v_4,v_5,v_6\}$ into the set of $P_4$ paths
$\mathcal C=\{(a_5,v_4,a_7,v_5)$, $(v_5,a_5,v_2,a_7)$, $(v_6,a_5,v_3,a_7)\}$.
Decompose the remaining edges between $\{a_9,a_{11}\}$ and $\{v_1,v_2,v_3,v_4,v_5,v_6\}$ into the set of $P_4$ paths
$\mathcal D=\{(a_9,v_3,a_{11},v_2)$, $(a_9,v_4,a_{11},v_5)$, $(a_{11},v_6,a_9,v_5)\}$.
Decompose the remaining edges between $\{a_2,a_4,a_6,a_8\}$ and $\{v_1,v_2,v_3,v_4,v_5,v_6\}$ along with the edges $(a_2,a_6)$ and $(a_4,a_1)$ into the set of $P_4$ paths
$\mathcal E=\{(v_4,a_2,a_6,v_5)$, $(v_5,a_2,v_2,a_6)$, $(a_1,a_4,v_6,a_8)$, $(v_1,a_4,v_2,a_8)$, $(v_1,a_8,v_3,a_4)$, $(v_3,a_6,v_4,a_8)$, $(v_3,a_2,v_6,a_6)\}$.
Decompose the remaining edges between $\{a_{10},a_{12}\}$ and $\{v_1,v_2,v_3,v_4,v_5,v_6\}$ into the set of $P_4$ paths
$\mathcal F=\{(v_4,a_{12},v_5,a_{10})$, $(a_{10},v_6,a_{12},v_1)$, $(v_1,a_{10},v_2,a_{12})\}$.

Let $P_1=P\setminus \{a_1,a_2,a_3,a_4,a_5,a_6,a_7,a_8,a_9,a_{10},a_{11},a_{12}\}$.

If $n\equiv 0,4$ (mod 6), let $m=\frac{|P_1|}{2}$. Partition $P_1$ into $m$ subsets $A_1,\ldots,A_{m}$, of size 2. For each $1\leq i\leq m$, decompose the edges between $\{v_1,v_2,v_3,v_4,v_5,v_6\}$ and $A_i$ into a set of $P_4$ paths $\mathcal A_i$ in a similar manner as Lemma~\ref{lem.d}.

If $n\equiv 1,3$ (mod 6), let $m=\lfloor \frac{|P_1|}{2} \rfloor$. Partition $P_1$ into $m$ subsets $A_1,\ldots,A_{m-1},A_{m}$, where $A_1,\ldots,A_{m-1}$ are of size 2 and $A_{m}$ is of size 3. For any $1\leq i\leq m-1$, decompose the edges between $\{v_1,v_2,v_3,v_4,v_5,v_6\}$ and $A_i$ into a set of $P_4$ paths $\mathcal A_i$ in a similar manner as Lemma~\ref{lem.d}. Decompose the edges between $\{v_1,v_2,v_3,v_4,v_5,v_6\}$ and $A_{m}$ into a set of $P_4$ paths $\mathcal A_{m}$ in a similar manner as Lemma~\ref{lem.e}.

Decompose the edges between the vertices in the set $\{v_1,v_2,v_3,v_4,v_5,v_6\}$ into a set of $P_4$ paths $\mathcal V$ in a similar manner as Lemma~\ref{lem0.1} so that the decomposition includes the two paths $(v_6,v_2,v_3,v_5)$ and $(v_2,v_4,v_1,v_3)$.

Let $\mathcal P'=(\mathcal P\setminus\{(a_6,a_2,a_3,a_5),(a_2,a_4,a_1,a_3)\}) \cup \mathcal A \cup \mathcal B \cup \mathcal C \cup \mathcal D \cup \mathcal E \cup \mathcal F \cup (\bigcup \limits_{i=1}^{m}\mathcal A_i )\cup \mathcal V$. Then $(P',\mathcal P')$ is a uniquely 2-chromatic $P_4(n+6)$ with colour classes $C'_w=C_w \cup \{v_1,v_3,v_5\} $ and $C'_b=C_b\cup \{v_2,v_4,v_6\}$.
\end{Proof}

We will use Lemma~\ref{cong1_2} iteratively, but to begin the process we need to have two non-critical $P_4$ paths of the prescribed form.
Such paths are not present in the $P_4 (109)$ constructed in Theorem~\ref{Thm-109},
so we introduce a construction that produces a system with two non-critical paths having the required form.

\begin{Lemma}\label{cong0}
Let $(P,\mathcal P)$ be a uniquely 2-chromatic $P_4$ system of order $n$  where $n\equiv 1$ (mod 6) with colour classes $C_w$ and $C_b$ such that $|C_w| \geq 10$, $|C_b| \geq 6$, and $|C_w| >\frac{1}{3} |P| + 5$. If there exist non-critical $P_4$ paths $(a_1,a_2,b_1,b_2)$ and $(b_2,b_3,a_2,b_4)$ in $\mathcal P$ such that $c(a_1)=c(a_2)=w$ and $c(b_1)=c(b_2)=c(b_3)=c(b_4)=b$, then there exists a uniquely 2-chromatic $P_4$ system of order $n+5$.
\end{Lemma}

\begin{Proof}
Let $\{a_1,a_2,a_3,a_4,a_5,a_6,a_7,a_8,a_9,a_{10}\} \subseteq C_w$ and $\{b_1,b_2,b_3,b_4,b_5,b_6\} \subseteq C_b$. Also, let $P'=P\cup \{v_1,v_2,v_3,v_4,v_5\}$. We will remove the two paths $(a_1,a_2,b_1,b_2)$ and $(b_2,b_3,a_2,b_4)$ from $\mathcal P$ so that we can reuse their edges in a different manner.

Let $\mathcal A= \{(a_2,a_1,v_1,a_3), (b_3,b_2,v_4,b_1), $ $(v_4,a_4,v_5,a_5),$ $ (v_2,b_3,v_1,b_4),$
$(a_2, v_4,a_6, v_3) \}$. Then $c(v_1)=c(v_3)=c(v_5)=b$ and $c(v_2)=c(v_4)=w$. Decompose the remaining edges between $\{a_1,a_3,a_5\}$ and $\{v_1,v_2,v_3,v_4,v_5\}$ into the set of $P_4$ paths $\mathcal B = \{(v_1,a_5,v_4,a_1)$, $(v_2,a_5,v_3,a_1)$, $(a_1,v_2,a_3,v_3)$, $(a_1,v_5,a_3,v_4)\}$. Decompose the remaining edges between $\{a_2,a_4,a_6\}$ and $\{v_1,v_2,v_3,v_4,v_5\}$ along with the edges $(a_2,b_4)$ and $(a_2,b_3)$ into the set of $P_4$ paths $\mathcal C=\{(b_3,a_2,v_1,a_4)$, $(b_4,a_2,v_2,a_6)$, $(v_1,a_6,v_5,a_2)$, $(v_2,a_4,v_3,a_2)\}$.

Decompose the remaining edges between $\{b_1,b_2,a_7\}$ and $\{v_1,v_2,v_3,v_4,v_5\}$ along with the edges $(b_2,b_1)$ and $(b_1,a_2)$ into the set of $P_4$ paths $\mathcal D=\{(a_2,b_1,b_2,v_1)$, $(v_2,b_1,v_3,a_7)$, $(a_7,v_1,b_1,v_5)$, $(v_3,b_2,v_5,a_7)$, $(b_2,v_2,a_7,v_4)\}$. Decompose the remaining edges between $\{b_3,b_4,a_8\}$ and $\{v_1,v_2,v_3,v_4,v_5\}$  into the set of $P_4$ paths
$\mathcal E=\{(v_1,a_8,v_5,b_4)$, $(b_3,v_3,a_8,v_4)$, $(b_4,v_4,b_3,v_5)$, $(a_8,v_2,b_4,v_3)\}$. Decompose the edges between $\{b_5,b_6,a_9\}$ and $\{v_1,v_2,v_3,v_4,v_5\}$ into the set of $P_4$ paths $\mathcal F=\{(v_1,a_9,v_5,b_6)$, $(b_5,v_3,a_9,v_4)$, $(v_5,b_5,v_4,b_6)$, $(a_9,v_2,b_6,v_3)$, $(b_6,v_1,b_5,v_2)\}$.

Let $P_1=P\setminus (\{a_1,a_2,a_3,a_4,a_5,a_6,a_7,a_8,a_9,a_{10},b_1,b_2,b_3,b_4,b_5,b_6\}$ and $m= \frac{|P_1|}{3} $. Partition $P_1$ into $m$ subsets $A_1,\ldots,A_{m}$ each of size 3 such that for any $1\leq i\leq m$, $A_i=\{a_1^i,a_2^i,a_3^i\}$ and $c(a_2^i)=w$ (this can be done because
$|C_w| >\frac{1}{3} |P| + 5$).

Decompose the edges between $\{v_1,v_2,v_3,v_4,v_5\}$ and $A_i$ into the set of $P_4$ paths
$\mathcal A_i=\{ (a_1^i,v_1,a_2^i,v_2)$, $(a_1^i,v_2,a_3^i,v_5)$, $(v_3,a_1^i,v_4,a_3^i)$, $(v_1,a_3^i,v_3,a_2^i)$, $(a_1^i,v_5,a_2^i,v_4)\}$. Rename $a_{10}$ to $v_6$ and decompose the edges between the vertices in the set $\{v_1,v_2,v_3,v_4,v_5,v_6\}$ into a set of $P_4$ paths $\mathcal V$ in a similar manner as Lemma~\ref{lem0.1} so that the decomposition includes the two paths $(v_6,v_2,v_3,v_5)$ and $(v_2,v_4,v_1,v_3)$.

Let $\mathcal P'=(\mathcal P\setminus\{(a_1,a_2,b_1,b_2),(b_2,b_3,a_2,b_4)\}) \cup \mathcal A \cup \mathcal B \cup \mathcal C \cup \mathcal D \cup \mathcal E \cup \mathcal F\cup (\bigcup \limits_{i=1}^{m}\mathcal A_i )\cup \mathcal V$. Then $(P',\mathcal P')$ is a uniquely 2-chromatic $P_4(n+5)$ with colour classes $C'_w=C_w \cup \{v_2,v_4\} $ and $C'_b=C_b\cup \{v_1,v_3,v_5\}$.
\end{Proof}

\begin{Corollary}
There exists a uniquely 2-chromatic $P_4$ system of order $n$ for all $n\equiv 0$ (mod 6), and $n > 109$.
\end{Corollary}

\begin{Proof}
By Theorem~\ref{Thm-109}, there exists a uniquely 2-chromatic $P_4$ system $(P,\mathcal P)$ of order 109 containing two non-critical paths  $(27',1,2',2)$ and $(25',1,6',2)$ in $\mathcal R_1$ such that $c(25')=c(1)=w$ and $c(2)=c(2')=c(6')=c(27')=b$. Therefore, by Lemma~\ref{cong0}, there exists a uniquely 2-chromatic $P_4$ system $(P',\mathcal P')$ of order 114. Moreover, this uniquely 2-chromatic $P_4$ system of order 114 contains two non-critical paths  $(v_6,v_2,v_3,v_5)$ and $(v_2,v_4,v_1,v_3)$ in $\mathcal V$ as in Lemma~\ref{cong0}
where $c(v_1)=c(v_3)=c(v_5)=w$ and $c(v_2)=c(v_4)=c(v_6)=b$. Therefore, by Lemma~\ref{cong1_2} there exists a uniquely 2-chromatic $P_4$ system of order 120. Iteratively, we construct a uniquely 2-chromatic $P_4$ system of order $n$, for all $n\equiv 0$ (mod 6) and $n\geq 120$, by repeatedly applying Lemma~\ref{cong1_2}.
\end{Proof}

\begin{Lemma}\label{cong1_1}
Let $(P,\mathcal P)$ be a uniquely 2-chromatic $P_4$ system of order $n$  where $n\equiv 1,3,4$ (mod 6) with colour classes $C_w$ and $C_b$ such that $|C_w| \geq 6$ and $|C_b| \geq 6$. If there exist  non-critical $P_4$ paths $(a_1,a_2,b_1,b_2)$ and $(b_2,b_3,a_2,b_4)$ in $\mathcal P$ such that $c(a_1)=c(a_2)=w$ and $c(b_1)=c(b_2)=c(b_3)=c(b_4)=b$, then there exists a uniquely 2-chromatic $P_4$ system of order $n+6$.
\end{Lemma}

\begin{Proof}
Let $\{a_1,a_2,a_3,a_4,a_5,a_6\} \subseteq C_w$ and $\{b_1,b_2,b_3,b_4,b_5,b_6\} \subseteq C_b$. Also, let $P'=P\cup \{v_1,v_2,v_3,v_4,v_5,v_6\}$. We will remove the two paths  $(a_1,a_2,b_1,b_2)$ and $(b_2,b_3,a_2,b_4)$ in $\mathcal P$ so that we can reuse their edges in a different manner.

Let $\mathcal A= \{(a_1,a_2,v_1,a_3)$, $(v_1,b_3,v_2,b_4)$, $(v_1,b_5,v_6,b_6)$,
$(b_3, b_2,v_4, b_1)$, $(v_4,a_1,v_5,a_4)$, $(v_4,a_5,v_3,a_6)\}$. Then $c(v_1)=c(v_3)=c(v_5)=b$ and $c(v_2)=c(v_4)=c(v_6)=w$. Decompose the remaining edges between $\{a_2,a_3\}$ and $\{v_1,v_2,v_3,v_4,v_5,v_6\}$ along with the edges $(a_2,b_3)$ and $(a_2,b_4)$ into the set of $P_4$ paths
$\mathcal B = \{(b_4,a_2,v_2,a_3)$, $(b_3,a_2,v_6,a_3)$, $(a_2,v_4,a_3,v_3)$, $(v_3,a_2,v_5,a_3)\}$.
Decompose the remaining edges between $\{a_1,a_4\}$ and $\{v_1,v_2,v_3,v_4,v_5,v_6\}$ into the set of $P_4$ paths
$\mathcal C=\{(v_1,a_1,v_2,a_4)$, $(v_1,a_4,v_6,a_1)$, $(v_4,a_4,v_3,a_1)\}$.
Decompose the remaining edges between $\{a_5,a_6\}$ and $\{v_1,v_2,v_3,v_4,v_5,v_6\}$ into the set of $P_4$ paths
$\mathcal D=\{(v_1,a_5,v_2,a_6)$, $(v_1,a_6,v_6,a_5)$, $(v_4,a_6,v_5,a_5)\}$.
Decompose the remaining edges between $\{b_1,b_2\}$ and $\{v_1,v_2,v_3,v_4,v_5,v_6\}$ along with the edges $(b_1,b_2)$ and $(b_1,a_2)$ into the set of $P_4$ paths $\mathcal E=\{(a_2,b_1,b_2,v_1)$, $(b_2,v_2,b_1,v_1)$, $(v_6,b_2, v_5, b_1)$, $(v_6,b_1,v_3,b_2)\}$.
Decompose the remaining edges between $\{b_3,b_4\}$ and $\{v_1,v_2,v_3,v_4,v_5,v_6\}$ into the set of $P_4$ paths
$\mathcal F=\{(v_1,b_4,v_6,b_3)$, $(b_4,v_3,b_3,v_4)$, $(v_4,b_4,v_5,b_3)\}$.
Decompose the remining edges between $\{b_5,b_6\}$ and $\{v_1,v_2,v_3,v_4,v_5,v_6\}$ into the set of $P_4$ paths
$\mathcal G=\{(v_1,b_6,v_2,b_5)$, $(b_6,v_3,b_5,v_4)$, $(b_5,v_5,b_6,v_4)\}$.

Let $P_1=P\setminus \{a_1,a_2,a_3,a_4,a_5,a_6,b_1,b_2,b_3,b_4,b_5,b_6\}$.

If $n\equiv 4$ (mod 6), let $m=\frac{|P_1|}{2}$. Partition $P_1$ into $m$ subsets $A_1,\ldots,A_{m}$, of size 2. For each $1\leq i\leq m$, decompose the edges between $\{v_1,v_2,v_3,v_4,v_5,v_6\}$ and $A_i$ into a set of $P_4$ paths $\mathcal A_i$ in a similar manner as Lemma~\ref{lem.d}.

 If $n\equiv 1,3$ (mod 6), let $m=\lfloor \frac{|P_1|}{2} \rfloor$. Partition $P_1$ into $m$ subsets $A_1,\ldots,A_{m-1},A_{m}$, where $A_1,\ldots,A_{m-1}$ are of size 2 and $A_{m}$ is of size 3. For any $1\leq i\leq m-1$, decompose the edges between $\{v_1,v_2,v_3,v_4,v_5,v_6\}$ and $A_i$ into a set of $P_4$ paths $\mathcal A_i$ in a similar manner as Lemma~\ref{lem.d}. Decompose the edges between $\{v_1,v_2,v_3,v_4,v_5,v_6\}$ and $A_{m}$ into a set of $P_4$ paths $\mathcal A_{m}$ in a similar manner as Lemma~\ref{lem.e}. 
 
 Decompose the edges between the vertices in the set $\{v_1,v_2,v_3,v_4,v_5,v_6\}$ into a set of $P_4$ paths $\mathcal V$ in a similar manner as Lemma~\ref{lem0.1} so that the decomposition includes the two paths $(v_6,v_2,v_3,v_5)$ and $(v_2,v_4,v_1,v_3)$.

Let $\mathcal P'=(\mathcal P\setminus \{(a_1,a_2,b_1,b_2),(b_2,b_3,a_2,b_4)\}) \cup \mathcal A \cup \mathcal B \cup \mathcal C \cup \mathcal D \cup \mathcal E \cup \mathcal F \cup \mathcal G \cup (\bigcup \limits_{i=1}^{m}\mathcal A_i )\cup \mathcal V$. Then $(P',\mathcal P')$ is a uniquely 2-chromatic $P_4(n+6)$ with colour classes $C'_w=C_w \cup \{v_2,v_4,v_6\} $ and $C'_b=C_b\cup \{v_1,v_3,v_5\}$.
\end{Proof}

\begin{Corollary}
There exists a uniquely 2-chromatic $P_4$ system of order $n$ for all $n\equiv 1$ (mod 6), and $n\geq 109$.
\end{Corollary}

\begin{Proof}
By Theorem~\ref{Thm-109}, there exists a uniquely 2-chromatic $P_4$ system $(P,\mathcal P)$ of order 109 containing the non-critical paths $(27',1,2',2)$ and $(25',1,6',2)$ in $\mathcal R_1$ such that $c(25')=c(1)=w$ and $c(2)=c(2')=c(6')=c(27')=b$.
Therefore, by Lemma~\ref{cong1_1}, there exists a uniquely 2-chromatic $P_4$ system $(P',\mathcal P')$ of order 115.
Moreover, this uniquely 2-chromatic $P_4$ system of order 115 contains two non-critical paths $(v_6,v_2,v_3,v_5)$ and $(v_2,v_4,v_1,v_3)$ in $\mathcal V$
where $c(v_1)=c(v_3)=c(v_5)=w$ and $c(v_2)=c(v_4)=c(v_6)=b$. Therefore, by Lemma~\ref{cong1_2}, there exists a uniquely 2-chromatic $P_4$ system of order 121. Iteratively, we construct a uniquely 2-chromatic $P_4$ system of order $n$, for all $n\equiv 1$ (mod 6) and $n\geq 121$, by repeatedly applying Lemma~\ref{cong1_2}.
\end{Proof}

\begin{Lemma}\label{cong3}
Let $(P,\mathcal P)$ be a uniquely 2-chromatic $P_4$ system of order $n$  where $n\equiv 1$ (mod 6) with colour classes $C_w$ and $C_b$ such that $|C_w| \geq 3$, $|C_b| \geq 4$. If there exists a non-critical $P_4$ path $(b_1,b_2,a_1,b_3)$ in $\mathcal P$ such that $c(a_1)=w$ and $c(b_1)=c(b_2)=c(b_3)=b$, then there exists a uniquely 2-chromatic $P_4$ system of order $n+2$.
\end{Lemma}

\begin{Proof}
Let $\{a_1,a_2,a_3\} \subseteq C_w$ and $\{b_1,b_2,b_3,b_4\} \subseteq C_b$. Also, let $P'=P\cup \{v_1,v_2\}$. We will remove the path $(b_1,b_2,a_1,b_3)$ in $\mathcal P$ so that we can use its edges in a different manner.

Let $\mathcal A = \{(b_2,b_1,v_1,b_3),(v_1,a_1,v_2,a_2)\}$. Then $c(v_1)=w$ and $c(v_2)=b$. Decompose the remaining edges between $\{v_1,v_2\}$ and $\{a_1,a_2,a_3,b_1,b_2,b_3,b_4\}$ along with the edges $(a_1,b_2)$, $(a_1,b_3)$, and $(v_1,v_2)$ into the set of $P_4$ paths $\mathcal B=\{(a_1,b_2,v_1,v_2)$, $(a_1,b_3,v_2,b_2)$, $(a_2,v_1,a_3,v_2)$, $(b_1,v_2,b_4,v_1)\}$.

Let $P_1=P\setminus \{a_1,a_2,a_3,b_1,b_2,b_3,b_4\}$ and $m= \frac{|P_1|}{3} $. Partition $P_1$ into $m$ subsets $A_1,\ldots,A_{m}$ of size 3. For any $1\leq i\leq m$, decompose the edges between $\{v_1,v_2\}$ and $A_i=\{a_1^i,a_2^i,a_3^i\}$ into the set of $P_4$ paths $\mathcal A_i=\{(a_1^i,v_1,a_2^i,v_2),(a_1^i,v_2,a_3^i,v_1)\}$.

Let $\mathcal P'=(\mathcal P\setminus \{(b_1,b_2,a_1,b_3)\}) \cup \mathcal A \cup \mathcal B \cup  (\bigcup \limits_{i=1}^{m}\mathcal A_i )$. Then $(P',\mathcal P')$ is a uniquely 2-chromatic $P_4(n+2)$ with colour classes $C'_w=C_w \cup \{v_1\} $ and $C'_b=C_b\cup \{v_2\}$.
\end{Proof}

\begin{Corollary}
There exists a uniquely 2-chromatic $P_4$ system of order $n$ for all $n\equiv 3$ (mod 6), and $n > 109$.
\end{Corollary}
\begin{Proof}
By Theorem~\ref{Thm-109}, there exists a uniquely 2-chromatic $P_4$ system $(P,\mathcal P)$ of order 109 containing the non-critical path $(27',1,2',2) \in \mathcal R_1$ such that $c(1)=w$ and $c(2)=c(2')=c(27')=b$. Therefore, by Lemma~\ref{cong3}, there exists a uniquely 2-chromatic $P_4$ system $(P',\mathcal P')$ of order 111.

Moreover, both $(P,\mathcal P)$ and $(P',\mathcal P')$ contain two non-critical paths  $(27',3,2',4)$ and $(25',3,6',4)$ in $\mathcal R_3$ such that $c(25')=c(3)=w$ and $c(6')=c(4)=c(2')=c(27')=b$. Therefore, by Lemma~\ref{cong1_1}, there exists a uniquely 2-chromatic $P_4$ system $(P'',\mathcal P'')$ of order 117.

This uniquely 2-chromatic $P_4$ system of order 117 contains two non-critical paths  $(v_6,v_2,v_3,v_5)$ and $(v_2,v_4,v_1,v_3)$
where $c(v_1)=c(v_3)=c(v_5)=w$ and $c(v_2)=c(v_4)=c(v_6)=b$. Therefore, by Lemma~\ref{cong1_2}, there exists a uniquely 2-chromatic $P_4$ system of order 123. Iteratively, we construct a uniquely 2-chromatic $P_4$ system of order $n$, for all $n\equiv 3$ (mod 6) and $n\geq 123$, by repeatedly applying Lemma~\ref{cong1_2}.
\end{Proof}

\begin{Lemma}\label{cong4}
Let $(P,\mathcal P)$ be a uniquely 2-chromatic $P_4$ system of order $n$, where $n\equiv 1$ (mod 6) with colour classes $C_w$ and $C_b$ such that $|C_w| \geq 3$, $|C_b| \geq 4$. If there exists a non-critical $P_4$ path $(b_1,b_2,a_1,b_3)$ in $\mathcal P$ such that $c(a_1)=w$ and $c(b_1)=c(b_2)=c(b_3)=b$, then there exists a uniquely 2-chromatic $P_4$ system of order $n+3$.
\end{Lemma}

\begin{Proof}
Let $\{a_1,a_2,a_3\} \subseteq C_w$ and $\{b_1,b_2,b_3\} \subseteq C_b$. Also, let $P'=P\cup \{v_1,v_2,v_3\}$. We will remove the path $(b_1,b_2,a_1,b_3)$ in $\mathcal P$ so that we can use its edges in a different manner.

Let $\mathcal A = \{(b_2,b_1,v_1,b_3),(v_1,a_1,v_2,a_2),(v_2,b_2,v_3,b_3)\}$. Then $c(v_1)=c(v_3)=w$ and $c(v_2)=b$. Decompose the remaining edges between $\{v_1,v_2,v_3\}$ and $\{a_1,a_2,a_3,b_1,b_2,b_3\}$ along with the edges $(a_1,b_2)$, $(a_1,b_3)$, and $(v_1,v_2),(v_1,v_3),(v_2,v_3)$ into the set of $P_4$ paths $\mathcal B=\{a_1,b_3,v_2,v_1)$, $(b_2,a_1,v_3,v_2)$, $(b_1,v_3,a_2,v_1)$, $(b_1,v_2,a_3,v_1)$, $(a_3,v_3,v_1,b_2)\}$.

Let $P_1=P\setminus \{a_1,a_2,a_3,b_1,b_2,b_3\}$ and $m=\lfloor \frac{|P_1|}{2} \rfloor$. Partition $P_1$ into $m$ subsets $A_1,\ldots,A_{m-1}$, $A_m$, where $A_1,\ldots,A_{m-1}$ are of size 2 and $A_m$ is of size 3 such that for any $1\leq i\leq m-1$, $A_i=\{a_1^i,a_2^i\}$, $A_m=\{a_1^m,a_2^m,a_3^m\}$, and $c(a_2^m)=b$.
For any $1\leq i\leq m-1$, decompose the edges between $\{v_1,v_2,v_3\}$ and $A_i=\{a_1^i,a_2^i\}$ into the set of $P_4$ paths $\mathcal A_i=\{(a_1^i,v_3,a_2^i,v_2)$, $(v_2,a_1^i,v_1,a_2^i)\}$. Decompose the edges between $\{v_1,v_2,v_3\}$ and $A_m$ into the set of $P_4$ paths $\mathcal A_m=\{(v_1,a_1^m,v_2,a_2^m)$, $(a_3^m,v_1,a_2^m,v_3)$, $(a_1^m,v_3,a_3^m,v_2)\}$.

Let $\mathcal P'=(\mathcal P\setminus \{(b_1,b_2,a_1,b_3)\}) \cup \mathcal A \cup \mathcal B \cup  (\bigcup \limits_{i=1}^{m}\mathcal A_i )$. Then $(P',\mathcal P')$ is a uniquely 2-chromatic $P_4(n+3)$ with colour classes $C'_w=C_w \cup \{v_1,v_3\} $ and $C'_b=C_b\cup \{v_2\}$.
\end{Proof}

\begin{Corollary}
There exists a uniquely 2-chromatic $P_4$ system of order $n$ for all $n\equiv 4$ (mod 6), and $n > 109$.
\end{Corollary}

\begin{Proof}
By Theorem~\ref{Thm-109}, there exists a uniquely 2-chromatic $P_4$ system $(P,\mathcal P)$ of order 109 containing the non-critical path $(27',1,2',2) \in \mathcal R_1$ such that $c(1)=w$ and $c(2)=c(2')=c(27')=b$. Therefore, by Lemma~\ref{cong4}, there exists a uniquely 2-chromatic $P_4$ system $(P',\mathcal P')$ of order 112.

Moreover, both $(P,\mathcal P)$ and $(P',\mathcal P')$ contain two non-critical paths $(27',3,2',4)$ and $(25',3,6',4)$ in $\mathcal R_3$ such that $c(25')=c(3)=w$ and $c(6')=c(4)=c(2')=c(27')=b$. Therefore, by Lemma~\ref{cong1_1}, there exists a uniquely 2-chromatic $P_4$ system $(P'',\mathcal P'')$ of order 118.

This uniquely 2-chromatic $P_4$ system of order 118 contains two non-critical paths  $(v_6,v_2,v_3,v_5)$ and $(v_2,v_4,v_1,v_3)$
where $c(v_1)=c(v_3)=c(v_5)=w$ and $c(v_2)=c(v_4)=c(v_6)=b$. Therefore, by Lemma~\ref{cong1_2}, there exists a uniquely 2-chromatic $P_4$ system of order 124. Iteratively, we construct a uniquely 2-chromatic $P_4$ system of order $n$, for all $n\equiv 4$ (mod 6) and $n\geq 124$, by repeatedly applying Lemma~\ref{cong1_2}.
\end{Proof}

Therefore, we have the following theorem.

\begin{Theorem}
There exists a uniquely 2-chromatic $P_4$ system of order $n$, for all $n\equiv 0,1$ (mod 3), and $n\geq 109$.
\end{Theorem}

\section{Acknowledgements}
D.~A.~Pike acknowledges research support from NSERC (grant number RGPIN-04456-2016).


\end{document}